\documentclass[smallextended,numbook,runningheads]{svjour3_hidden}     
\smartqed  
\usepackage[margin=3cm]{geometry} 

\usepackage{graphicx}
\usepackage{amsmath}
\usepackage{epstopdf} 
\usepackage{mathptmx}      
\journalname{}

\usepackage{stmaryrd}
\usepackage{psfrag}
\usepackage{subfig}
\usepackage{graphicx}
\usepackage{etoolbox}
\usepackage{amssymb}
\usepackage[show]{ed}
\usepackage[english]{babel}
\usepackage{tikz}
\usepackage{appendix}
\usepackage{mathtools} 
\usepackage{bbm}
\usepackage{enumerate}
\usepackage{hyperref}
\spnewtheorem{assumption}{Assumption}[section]{\bf}{\it}
\usepackage[numbers]{natbib}
\mathtoolsset{showonlyrefs} 
\newcommand{\pdiv}[2]{\frac{\partial #1}{\partial #2}}
\newcommand{\pdivl}[2]{{\partial #1}/{\partial #2}}
\newcommand{\imag}{\mathrm{i}}

\newcommand{\dd}{\ensuremath{\operatorname{d}\!}}
\newcommand*{\e}{\mathrm{e}}

\newcommand{\bfn}{\vec{n}}
\newcommand{\bfd}{\vec{d}}
\newcommand{\bfx}{\vec{x}}
\newcommand{\bfy}{\vec{y}}

\newcommand{\R}{\mathbb{R}}
\newcommand{\C}{\mathbb{C}}

\newcommand{\Z}{\mathbb{Z}}

\allowdisplaybreaks[4]

\newcommand{\McansChSix}{{M}}
\newcommand{\MTayChSix}{{N_T}}
\newcommand{\MtruncChSeven}{{\hat{N}}}
\newcommand{\Nquad}{{N_Q}}
\newcommand{\NsidesChSix}{{N}}
\newcommand{\uinc}{u^{i}}
\newcommand{\uiPW}{u^{i}_{\alpha}}
\newcommand{\gHerg}{g_{\mathop{\rm Herg}}}
\newcommand{\uiHerg}{u^i_{\mathop{\rm Herg}}}

\newcommand{\boundary}{{\partial\Omega}}

\begin{document}
	
\newtoggle{imp}
\toggletrue{imp}

\title{Numerically stable computation of embedding formulae for scattering by polygons\thanks{This work was supported by an EPSRC PhD studentship and an Australian Bicentennial Scholarship, both awarded to Andrew Gibbs.
}}

\author{A. Gibbs         \and
        S. Langdon			\and
        A. Moiola
}

\institute{A. Gibbs \at
              Department of Computer Science, 
              K.U. Leuven, 
              Celestijnenlaan 200A, 
              BE-3001 Leuven, 
              Belgium  \\
              \email{andrew.gibbs@cs.kuleuven.be}         
           \and
           S. Langdon \at
              Department of Mathematics and Statistics,
              University of Reading,
              Whiteknights, PO Box 220,
              Berkshire, RG6 6AX, UK\\
              Tel.: +44 118 378 5021\\
              \email{s.langdon@reading.ac.uk}
              \and
              A. Moiola \at
              Dipartimento di Matematica ``F. Casorati'',
              Universit\`a degli studi di Pavia,
              Via Ferrata 5, 27100 Pavia, Italy\\
              Tel.: +39 0382 985656\\
              \email{andrea.moiola@unipv.it}
}

\date{ }

\maketitle

\begin{abstract}
For problems of time-harmonic scattering by polygonal obstacles, embedding formulae provide a useful means of computing the far-field coefficient induced by any incident plane wave, given the far-field coefficient of a relatively small set of canonical problems. The number of such problems to be solved depends only on the geometry of the scatterer. Whilst the formulae themselves are exact in theory, any implementation will inherit numerical error from the method used to solve the canonical problems. This error can lead to numerical instabilities. Here, we present an {effective} approach to identify and regulate these instabilities. 
This approach is subsequently extended to {the case where the incident wave is a} Herglotz wave function, and we suggest how this could potentially remove frequency dependence of a T-matrix method.

\keywords{Embedding formulae \and Far-field  \and Helmholtz {equation} \and {Scattering}}
\subclass{35J05{,
74J20, 
76Q05, 
78A45}} 
\end{abstract}

\section{Introduction}\label{s:Intro}

In problems of {two-dimensional} time-harmonic acoustic scattering, obtaining a full characterisation of the scattering properties of an obstacle may require a representation of the far-field behaviour induced by a large set (perhaps even a continuum) of incident plane waves, each of the form
\begin{equation}\label{def:PWinc}
\uiPW(\bfx):=\e^{-\imag k(x_1\cos\alpha+x_2\sin\alpha)},\quad\bfx {=(x_1,x_2)}\in \R^2,
\end{equation}
with incident angle $\alpha\in[0,2\pi)$ and wavenumber $k>0$.
In this paper we consider {solutions to the Helmholtz equation}
\begin{equation}\label{Helm}
(\Delta+k^2)u_\alpha=0\quad\text{in }\R^2\setminus{\overline{\Omega}},
\end{equation}
together with the Dirichlet boundary condition
\begin{equation}\label{Dir}
u_\alpha=0\quad\text{on }\boundary,
\end{equation}
where {$\Omega$ is a simply connected{, bounded,} open set in $\R^2$, corresponding to the interior of a convex polygonal scatterer}{.} 
{Here $u_\alpha=\uiPW+u_\alpha^s$ denotes the total field, where $\uiPW$ is a (given) incident plane wave as in} (\ref{def:PWinc}) (combinations of these will be discussed in~\S\ref{sec:HergEmbed}),
and the scattered field $u_\alpha^s$ satisfies the Sommerfeld radiation condition (see, e.g., \cite[Definition~2.4]{CoKr:13})
\begin{equation}\label{SRC}
\pdiv{u^s_{\alpha}}{r}(\bfx)-\imag k u^s_{\alpha}(\bfx) = o(r^{-1/2}), \quad\text{as }r\rightarrow\infty,
\end{equation}
where $r:=|\bfx|$.

The \emph{far-field diffraction coefficient} (also called the \emph{far-field pattern} (e.g., \cite{CoKr:13}), \emph{far-field directivity} (e.g., \cite{ShCr:10}), or simply \emph{far-field coefficient} (e.g., \cite{Bi:16}))
is of practical interest, and is central to this paper.  
Intuitively, it describes the distribution of energy of the scattered field, measured far away from the scatterer.
Denoting the far-field coefficient at observation angle $\theta$ (where $\bfx = r(\cos\theta,\sin\theta)$) by $D(\theta,\alpha)$, it is formally defined by (see, e.g., \cite[Theorem~2.6]{CoKr:13})
\begin{equation}\label{FFlabelTest}
u_\alpha^s(\bfx)= \frac{\e^{\imag(kr +\pi/4)}}{\sqrt{2\pi k r}}\bigg(D(\theta,\alpha)+\mathcal{O}(r^{-1})\bigg)\quad \text{for}\quad r\rightarrow\infty.
\end{equation}
In what follows, we will make use of the following representation of the far-field coefficient (see, e.g., \cite[(2.14)]{CoKr:13}),
\begin{equation}\label{FFrep}
D(\theta,\alpha)=-\int_\boundary \e^{-\imag k(y_1\cos\theta+ y_2\sin\theta)}\pdiv{u_\alpha}{\bfn}(\bfy)\dd{s}(\bfy),
\end{equation}
where $\bfy=(y_1,y_2)\in\R^2$, and the normal derivative (Neumann data) is 
\[
\pdiv{u_\alpha}{\bfn}(\bfx):=\bfn(\bfx)\cdot\nabla_{\bfx} u_\alpha(\bfx),
\]
where $\bfn(\bfx)$ denotes the outward normal to $\boundary$ at $\bfx$. 
In practice, it is typical to solve \eqref{Helm}--\eqref{SRC} using a numerical solver, for example a \emph{Boundary Element Method}, 
from which the Neumann data can be obtained and inserted into \eqref{FFrep} to obtain the far-field coefficient.  In this paper, we are primarily interested in efficient and stable 
numerical approximation of $D(\theta,\alpha)$ for many pairs 
$(\theta,\alpha)\in[0,2\pi)^2$, for the case where $\Omega$ is a convex rational polygon (see Definition~\ref{def:RationalPolygons}), using \emph{embedding formulae}.

Embedding formulae have been developed for a range of geometries in recent years. The starting point is always the same, to seek an operator which:
\begin{enumerate}[(i)]
\item commutes with the Helmholtz equation~\eqref{Helm}, i.e.,\ maps solutions of the Helmholtz equation to solutions of the Helmholtz equation;
\item preserves boundary conditions (in our case, the sound-soft condition~\eqref{Dir});
\item maps plane waves to zero, hence mapping the total field to the scattered field.
\end{enumerate}
In general this will lead to radiating solutions of \eqref{Helm}-\eqref{Dir} which contain \emph{overly singular} terms, violating the so-called \emph{Meixner condition}, which ensures that no energy radiates from the corners of the scattering obstacle, and is required for uniqueness (see, e.g., \cite{Bi:06} or \cite{CrSh:05} for definitions and further details). Embedding formulae typically follow by choosing a linear combination of these overly singular terms {which sum} to zero, hence constructing a solution which satisfies the Meixner condition{. B}y a uniqueness argument{, this} must be the solution to the original scattering problem (see, e.g., \cite{Bi:06} or \cite{CrSh:05}).

The power of {embedding formulae} is that they enable us to explicitly represent the far-field coefficient $D(\cdot,\alpha)$ for any incident angle $\alpha\in[0,2\pi)$ in terms of the far-field coefficients 
$D(\cdot,\alpha_m)$ for just a few incident angles $\alpha_m$, $m=1,\ldots,M$, where $M$ depends only on the geometry of the polygon~\cite{Bi:06}.  In principle, this thus enables one to obtain a full characterisation of the scattering properties of the obstacle through solving only a small number of scattering problems; since numerical methods are usually required to approximate far-field coefficients, these formulae can potentially lead to a significant reduction in computational cost.

However, whilst the formulae themselves are in principle exact, any implementation will inherit numerical error from the method used to approximate the solutions of the canonical problems {(i.e.\ those for the incident angles $\alpha_1,\ldots,\alpha_M$)}.
We demonstrate below that even a small error in the numerical method used to construct the solutions to the canonical problems can lead to arbitrarily large error at certain observation points in the approximate far-field pattern, if obtained by \emph{naively} inserting the approximate canonical solutions into the embedding formulae.  We provide upper and lower bounds on the amplification of this error, in terms of the error in the solver used to obtain the canonical solutions.  We then suggest a new formulation, based on a series of carefully chosen Taylor expansions, which we can show to be stable, a claim {supported} by numerical examples. 

As mentioned above, in this paper we are primarily concerned with the practical implementation of embedding formulae for two-dimensional rational polygons (see, e.g., \cite{Bi:06} or~\cite{CrSh:05}).  The formulae developed in~\cite{CrSh:05}, extended to three-dimensional polyhedra in \cite{SkCrShVa:10}, contain $M$ edge Green functions, but due to the singular nature of the Green's functions, even the two-dimensional version of this approach is challenging to implement.  The approach of~\cite{Bi:06} instead requires the far-field coefficient induced by $M$ distinct incident plane waves, as outlined above.
The issues that arise with numerical implementation of the embedding formulae occur in both the Green's function approach of \cite{CrSh:05} (on top of the issues associated with computing the highly singular edge Green's functions) and the alternative approach of \cite{Bi:06}, although here we focus on the ideas from \cite{Bi:06}.
Embedding formulae also exist for a range of infinite and semi-infinite structures in two dimensions, such as infinite rational wedges and the semi-infinite strip \cite{CrSh:05}, and continue to generate interest in the academic community; problems of polygonal structures inside waveguides were considered in~\cite{Bi:16} and problems on a family of curved obstacles were considered in~\cite{MoBiCh:16}.

Given that embedding formulae have typically been used to provide far-field representations for problems of plane wave incidence, it is a natural extension to consider problems {where the incident field is a Herglotz wave function} (see Definition~\ref{hergs}), given that a Herglotz wave function is a weighted integral of plane waves. We will demonstrate that a naive implementation of \cite{Bi:06} also breaks down {in such a case}, due to numerical instabilities at certain points in the integral. Using our alternative numerically stable approach, we are able to overcome this problem. Moreover, we demonstrate that the Herglotz wave function incidence is applicable to problems of practical interest, in particular as an enhancement to \emph{Tmatrom} (see \cite{GaHa:16}), a recently developed numerically stable T-matrix software package.

An outline of this paper is as follows: in~\S\ref{problem}, we introduce embedding formulae, and the problems to which they apply. We also explain why and where numerical instabilities can occur when the formulae are implemented in practice. In~\S\ref{stability} we derive an alternative{, numerically stable,} formulation of the embedding formulae, and in~\S\ref{sec:imp} key implementational aspects of this formulation are discussed. Finally, in \S\ref{sec:HergEmbed} we extend the embedding formulae to Herglotz wave function incidence, and present an application to Tmatrom.

\section{Embedding formulae}\label{problem}
In this paper we assume that $\Omega$ is a convex and \emph{rational} $ \NsidesChSix $-sided polygon $\Omega$, which we define here, alongside related parameters for the embedding formulae that follow.
\begin{definition}[Rational polygons]\label{def:RationalPolygons}
An angle $\omega\in[0,2\pi)$ is \emph{rational} if it can be expressed as $\pi$ multiplied by a rational number. We describe an $\NsidesChSix$-sided polygon $\Omega$ as rational if its external angles 
$\{\omega_j\}_{j=1}^{\NsidesChSix}$ are all rational angles.
	For a rational polygon $\Omega$, we denote by $p$ the smallest {positive} integer such that 
	\begin{equation}\label{def:pQIrreg}
	\frac{\pi}{p}\text{ divides }\omega_j\text{ exactly },\quad\text{for }j=1,\ldots,\NsidesChSix,
	\end{equation}
	whilst $\{q_j\}_{j=1}^{\NsidesChSix}$ denotes the set of integers such that
	\begin{equation}\label{def:qQIrreg}
	\frac{q_j\pi}{p}=\omega_j\text{ for }j=1,\ldots,\NsidesChSix.
	\end{equation}
	We say a polygon is \emph{quasi-regular} if each of its external angles are equal.  In such a case 
	we define $q:=q_1=\ldots=q_N$.
\end{definition}
We remark that whilst embedding formulae for two-dimensional infinite wedges with angles which are an irrational multiple of $\pi$ were developed in~\cite{ShCr:10}, to the best knowledge of the authors no such formulae yet exist for polygons. Hence in this paper we restrict our attention to polygons for which each angle is a rational multiple of $\pi$. Naturally, any irrational angle may be approximated arbitrarily well by a rational angle, although the analysis of such an approximation is outside of the scope of this paper. {We choose our coordinate system so that} at least one side of $\Omega$ is parallel to the x-axis, with the origin lying inside~$\Omega$.

As outlined above, the embedding formulae of \cite[(3.4)]{Bi:06} enable us to (in principle) compute the far-field coefficient $D(\theta,\alpha)$ for a range of $\alpha\in[0,2\pi)$, given the solution for just a few canonical problems corresponding to incident angles $\alpha_m$, $m=1,\ldots,\McansChSix$.
The parameter $\McansChSix$ depends only on the geometry of $\Omega$. 
In the case of a quasi-regular polygon, it follows from \cite[\S3]{Bi:06} that
\[
\McansChSix=\NsidesChSix(q-1),
\]
which corresponds to $q-1$ incident waves for each corner of $\Omega$. The choice of $q$ is discussed further in Remark~\ref{Mchoice} below. In the case of an irregular rational polygon, it follows from \cite[\S3.2]{Bi:06} that
\[
\McansChSix=\sum_{j=1}^\NsidesChSix(q_j-1),
\]
which corresponds to $q_j-1$ incident waves for the $j$th corner of $\Omega$. Exactly which incident angles to choose for the $\McansChSix$ canonical problems is discussed in Remark~\ref{Mchoice2}, whilst further constraints are introduced in~\S\ref{MV1}.

\begin{remark}\label{rem:greekLetters}
Throughout the paper, {{the symbols} $\alpha$ and $\theta$} are to be interpreted as 
{angles representing points} on the unit circle, so for example $\theta$ {is identified} with $\theta+2\pi${. W}e say that $\theta_1$ is \emph{close to} $\theta_2$ if $|\e^{\imag\theta_1}-\e^{\imag\theta_2}|<\epsilon$ for some $\epsilon>0$ small.
\end{remark}

In what follows we will require the following result (see, e.g., \cite[Theorem~3.15]{CoKr:13}): 
\begin{theorem}[Reciprocity relation]\label{TH:recip}
{For $\alpha,\theta\in[0,2\pi)$ and} a sound-soft obstacle $\Omega$, we have the identity:
\[D(\theta,\alpha)=D(\alpha,\theta).\]
\end{theorem}
We solve \eqref{Helm}--\eqref{SRC} and compute \eqref{FFrep} for canonical incident angles 
\[ \{\alpha_1,\ldots,\alpha_\McansChSix\}=:A_\McansChSix. \]
We define
\begin{equation}\label{def:Lambda}
\Lambda(\theta,\alpha):=\cos(p\theta)-(-1)^p\cos(p\alpha)
\end{equation}
and, recalling Theorem~\ref{TH:recip},
\begin{equation}\label{def:Dhat}
\hat{D}(\theta,\alpha):=\Lambda(\theta,\alpha){D}(\theta,\alpha) = (-1)^{p+1}\hat{D}(\alpha,\theta).
\end{equation}
The embedding formula follows from \cite[(3.4)]{Bi:06}:
\begin{equation}\label{emb0}
{D}(\theta,\alpha)=\frac{\sum_{m=1}^\McansChSix B_m(\alpha)\Lambda(\theta,\alpha_m){D}(\theta,\alpha_m)}{\Lambda(\theta,\alpha)}, \quad\text{for }(\theta,\alpha)\in[0,\pi)^2,
\end{equation}
where $[B_m]_{m=1}^\McansChSix\in\C^\McansChSix$ solves the system of equations
\begin{equation}\label{Bm_def}
\sum_{m=1}^\McansChSix B_m(\alpha)\hat{D}(\alpha_n,\alpha_m)=(-1)^{p+1}\hat{D}(\alpha,\alpha_n),
\end{equation}
for $n=1,\ldots,\McansChSix$. 
For $(\theta,\alpha)\in[0,2\pi)^2$ such that $\Lambda(\theta,\alpha)\neq0$, the representation \eqref{emb0} can, in principle, be evaluated explicitly to obtain the far-field coefficient $D(\theta,\alpha)$. {As explained in \cite{Bi:06}, f}or the case $\Lambda(\theta,\alpha)=0$, L'Hopit\^{a}l's rule may be applied to obtain $D(\theta,\alpha)$, with a second application in the sub-case $[\pdivl{\Lambda}{\theta}](\theta,\alpha)=0$. Put formally, (at least) one application of L'Hopit\^{a}l's rule is required for $\theta$ in the set
\begin{align}\label{LHoppyLopps}
\Theta_\alpha&:=\{\theta\in[0,2\pi):\Lambda(\theta,\alpha)=0\}\nonumber\\
&=\left\{\begin{array}{lll}
\{\theta\in[0,2\pi): \theta=\pm\alpha+({2n+1})\pi/{p},\quad& n\in\mathbb{Z}\},&p\text{ odd},\\
\{\theta\in[0,2\pi): \theta=\pm\alpha+{2n}\pi/{p},\quad& n\in\mathbb{Z}\},&p\text{ even,}
\end{array}\right.
\end{align}
with a second application if $[\pdivl{\Lambda}{\theta}](\theta,\alpha)=0$, i.e.\ if 
\begin{equation}\label{LHoppyLopps2}
{\theta\in \Theta_*:=\{\theta\in[0,2\pi):[\pdivl{\Lambda}{\theta}](\theta,\alpha)=0\}=\{\theta\in[0,2\pi):\theta=n\pi/p,\; n\in \mathbb Z\}}
\end{equation}
also. Considering Remark \ref{rem:greekLetters}, the set $\Theta_\alpha$ 
contains $2p$ elements, unless $\alpha=(2n+1)\pi/p$, $n\in\mathbb{Z}$, for $p$ odd, or $\alpha=2\pi n/p$, $n\in\mathbb{Z}$, for $p$ even; in either of these cases, 
$\Theta_\alpha$ contains $p$ elements (in which case these elements are also in $\Theta_*$).  We note that $\Theta_{\alpha}\subset\Theta_*$ if and only if $\alpha\in\Theta_*$, though this does not imply that $\Lambda(\theta,\alpha)=0$ for all $(\theta,\alpha)\in \Theta_*\times \Theta_*$.

We now summarise the process of implementing embedding formulae.
\begin{enumerate}[(i)]
\item Given an obstacle $\Omega$, compute $p$, $\{q_j\}_{j=1}^{ \NsidesChSix }$ and $\McansChSix$ (details in Remark~\ref{Mchoice} and Remark~\ref{Mchoice2}  below).
\item {Solve the Helmholtz boundary value problem \eqref{Helm}--\eqref{SRC} and c}ompute the far-field coefficient $D(\cdot,\alpha_m)$ for $\McansChSix$ distinct incident angles~$\alpha_m$, $m=1,\ldots,M$.
\item Given an incident angle $\alpha$, solve the $\McansChSix\times\McansChSix$ system \eqref{Bm_def} to determine the coefficients~$B_m(\alpha)$.
\item Obtain $D(\theta,\alpha)$ from \eqref{emb0}, using \eqref{LHoppyLopps} and \eqref{LHoppyLopps2} to determine if application(s) of L'Hopit\^{a}l's rule is/are required.
\end{enumerate}
Note that for a different incident angle $\alpha$, steps~(i) and~(ii) need not be repeated.

\begin{remark}[Reduction in number of canonical solves required]\label{Mchoice}
	In the case of quasi-regular polygons 
	we have the elementary formula
	\begin{equation}\label{QRangles}
	\omega_j=\pi\frac{\NsidesChSix+2}{\NsidesChSix}, \quad\text{for }j=1,\ldots,\NsidesChSix,
	\end{equation}
	from which the parameters $p$ and $q$ can be determined immediately. In \cite[\S3]{Bi:06} the number of distinct incident angles required for quasi-regular polygons is stated as $\McansChSix=\NsidesChSix(\NsidesChSix+1)$. This is based on the idea that $q~-~1$ canonical solutions are needed for each corner of the polygon, hence if we take the {simplest} choice (considering \eqref{QRangles}) of $q=\NsidesChSix+2$, this yields
	\begin{equation}\label{OldMdef}
	\McansChSix=\NsidesChSix(q-1)=\NsidesChSix(\NsidesChSix+1).
	\end{equation}
	However, the {simplest}
	choice of $q=\NsidesChSix+2$ may not be the best choice. Noting \eqref{QRangles} for even $\NsidesChSix$, we may instead choose $q=(\NsidesChSix+2)/2$ and $p=\NsidesChSix/2$ and thus \eqref{OldMdef} becomes
	\[
	\McansChSix=\NsidesChSix(q-1)=\NsidesChSix\left(\frac{\NsidesChSix+2}{2}-1\right)=\frac{\NsidesChSix^2}{2}.
	\]
	For even $\NsidesChSix$, this reduces the number of solves, compared to the choice of~\cite[(3.3)]{Bi:06}, by more than half.
	Although the choice suggested by \cite{Bi:06} is still valid, the system \eqref{emb0} may become redundant when implemented numerically, if $\McansChSix$ is larger than is necessary. 
\end{remark}

\begin{remark}\label{Mchoice2}
	In \cite{Bi:06}, no constraint is placed on the incident angles $A_{\McansChSix}$, other than that they are distinct. Numerical experiments have shown that close clustering of these incident angles may cause the embedding formulae to become inaccurate. Given that the embedding formulae for polygons extend naturally from that of a rational wedge, one would expect that if at least $q-1$ incident angles can `{see}' any given corner of $\Omega$, that is $\bfn(\bfx)\cdot(\cos(\alpha),\sin(\alpha))>0$, then the embedding formulae would hold. Numerical experiments suggest that this is sufficient, and we will assume that this constraint holds for the remainder of the paper.
	\end{remark}

We remind the reader that this paper is primarily concerned with efficient and stable approximation of the far-field coefficient $D(\theta,\alpha)$, for any $(\theta,\alpha)\in[0,2\pi)^2$. We now motivate this by demonstrating cases when a naive implementation of \eqref{emb0} can become unstable.  

Any numerical approximation to $u_\alpha$ will not be exact, and this can cause significant numerical instabilities close to $\Theta_\alpha$. To see why, we denote by 
$\mathcal{P}{D}(\theta,\alpha)$ some numerical approximation to ${D}(\theta,\alpha)$, 
and adopt the notational convention $\mathcal{P}\hat{D}:=\Lambda \mathcal{P}D$.  We can then define the embedding formula coupled with our numerical solver as
\begin{equation}\label{embP0}
D(\theta,\alpha)\approx\mathcal{E}_{\mathcal{P}}D(\theta,\alpha):=
\frac{\sum_{m=1}^\McansChSix b_m(\alpha)\mathcal{P}\hat{D}(\theta,\alpha_m)}{\Lambda(\theta,\alpha)},
\end{equation}
where $[b_m]_{m=1}^\McansChSix\approx [B_m]_{m=1}^\McansChSix$ solves the system of equations
\begin{equation}\label{approxBm}
\sum_{m=1}^\McansChSix b_m(\alpha){\mathcal{P}}\hat{D}(\alpha_n,\alpha_m)=(-1)^{p+1}{\mathcal{P}}\hat{D}(\alpha,\alpha_n),\quad\text{for }n=1,\ldots,\McansChSix.
\end{equation}
Successful implementation of the embedding formulae depends on the following assumption.
\begin{assumption}\label{b_m_unique}
\begin{enumerate}[(i)]
\item 
The {$M\times M$} system \eqref{Bm_def} is uniquely solvable for any $\alpha\in[0,2\pi)$.
\item 
The {$M\times M$} system \eqref{approxBm}  is uniquely solvable for any $\alpha\in[0,2\pi)$.
\end{enumerate}
\end{assumption}
To the best knowledge of the authors, the questions of precisely when Assumptions \ref{b_m_unique}(i) and~\ref{b_m_unique}(ii) hold are (both) still open.  {However, we have not found a case for which either part of Assumption \ref{b_m_unique} fails, moreover Assumption \ref{b_m_unique}(ii) has been observed to hold for a large range of numerical experiments.}

Throughout this paper, for the far-field solver $\mathcal{P}$, a modified version of MPSpack of \cite{BaBe:10} has been used. This version of MPSpack uses the method of fundamental solutions \cite{BaBe:08}; it was modified for Tmatrom \cite{GaHa:16}, and we have made a subsequent modification to that version allowing for fast computation of derivatives of the far-field coefficient (details given in given in \S\ref{sec:fastFF}). This will be required for our stable approximation in \S\ref{stability}. We note that the choice of $\mathcal{P}$ makes no difference to the points at which the embedding formulae become numerically unstable. 

It then follows that the error in the embedding formula is
\begin{equation}
\left|{D}(\theta,\alpha)-\mathcal{E}_{\mathcal{P}}D(\theta,\alpha)\right|
=\frac{1}{|\Lambda(\theta,\alpha)|}\left|\sum_{m=1}^\McansChSix\left[B_m(\alpha)\hat{D}(\theta,\alpha_m)-b_m(\alpha)\mathcal{P}\hat{D}(\theta,\alpha_m)\right]\right|.
\label{instab}
\end{equation}
It is clear from \eqref{instab} that for any amount of numerical error arising from the solver $\mathcal{P}$, a naive implementation of the embedding formula can lead to arbitrarily large error at points $(\theta,\alpha)$ such that $\Lambda(\theta,\alpha)\approx0$. We remark that in the alternative approach of \cite{CrSh:05}, which uses edge Green's functions, the division by $\Lambda(\theta,\alpha)$ still occurs (in particular see \cite[(4.4)]{CrSh:05}), hence the same instabilities will arise in practice.
The following lemma bounds the scaling of this error in terms of known parameters.

\begin{lemma}\label{HatBound}
	Given $p\in\mathbb{N}$,
	the following bounds hold for all $(\theta,\alpha)\in[0,2\pi)^2$:
	\begin{equation}\label{eq:HatBound}
	\frac{p^2}{8}|\theta-\theta_0||\theta_0-\theta_*|\leq|\Lambda(\theta,\alpha)|\leq{p^2}|\theta-\theta_0|\left(\frac{1}{2}|\theta-\theta_0|+|\theta_0-\theta_*|\right),
	\end{equation}
	where $\theta_0$ is the element of $\Theta_\alpha$ closest to $\theta$, and $\theta_*$ is the element of $\Theta_*$ closest to $\theta_0$ (\emph{close} in the sense of Remark \ref{rem:greekLetters}).
\end{lemma}
A simple interpretation of the lemma is that $1/\Lambda(\theta,\alpha)\sim|\theta-\theta_0|^{-1}$ when $\theta_0{\in\Theta_\alpha}$ is far from $\theta_*$, but $1/\Lambda(\theta,\alpha)\sim|\theta-\theta_0|^{-2}$ when $|\theta_0-\theta_*|$ and $|\theta-\theta_0|$ are small.
\begin{proof}
	Firstly, by the definition  \eqref{LHoppyLopps} of $\Theta_\alpha$ we have
	\begin{equation}\label{alpha_replace_theta0}
	\Lambda(\theta,\alpha)=
	{\Lambda(\theta,\alpha)-\Lambda(\theta_0,\alpha)}=
	\cos(p\theta)-\cos(p\theta_0),\quad\text{for }\theta_0\in\Theta_\alpha,
	\end{equation}
	and from standard trigonometric identities it follows that
	\begin{equation}\label{Tbound1}
	\left|{{\Lambda}{(\theta,\alpha)}}\right|
	=
	{2\left|\sin\left({p}(\theta-\theta_0)/{2}\right)\right|\cdot\left|\sin\left({p}(\theta+\theta_0)/{2}\right)\right|}.
	\end{equation}
	We first focus on the lower bound of \eqref{eq:HatBound}, for which we will require the inequality
	\begin{equation}\label{eq:dodgyInequality}
	|\sin(p(\theta_*+\theta_0)/2)|\leq|\sin(p(\theta+\theta_0)/2)|.
	\end{equation}
	\begin{figure}[hbpt]
		\centering
		{\psfrag{[x]}{$\theta$}\psfrag{[y]}[Bc][Bc]{}
			\psfrag{[x1]}[Bc][Bc]{$0$}\psfrag{[x2]}[Bc][Bc]{$\pi/2$}\psfrag{[x3]}[Bc][Bc]{$\pi$}\psfrag{[x4]}[Bc][Bc]{$3\pi/2$}\psfrag{[x5]}[Bc][Bc]{$2\pi$}\psfrag{[t]}[Bc][Bc]{$\alpha=0.3,$ $p=2$}\psfrag{[a]}[Bc][Bc]{$\Lambda(\theta,\alpha)$}
			\includegraphics[width=.495\textwidth]{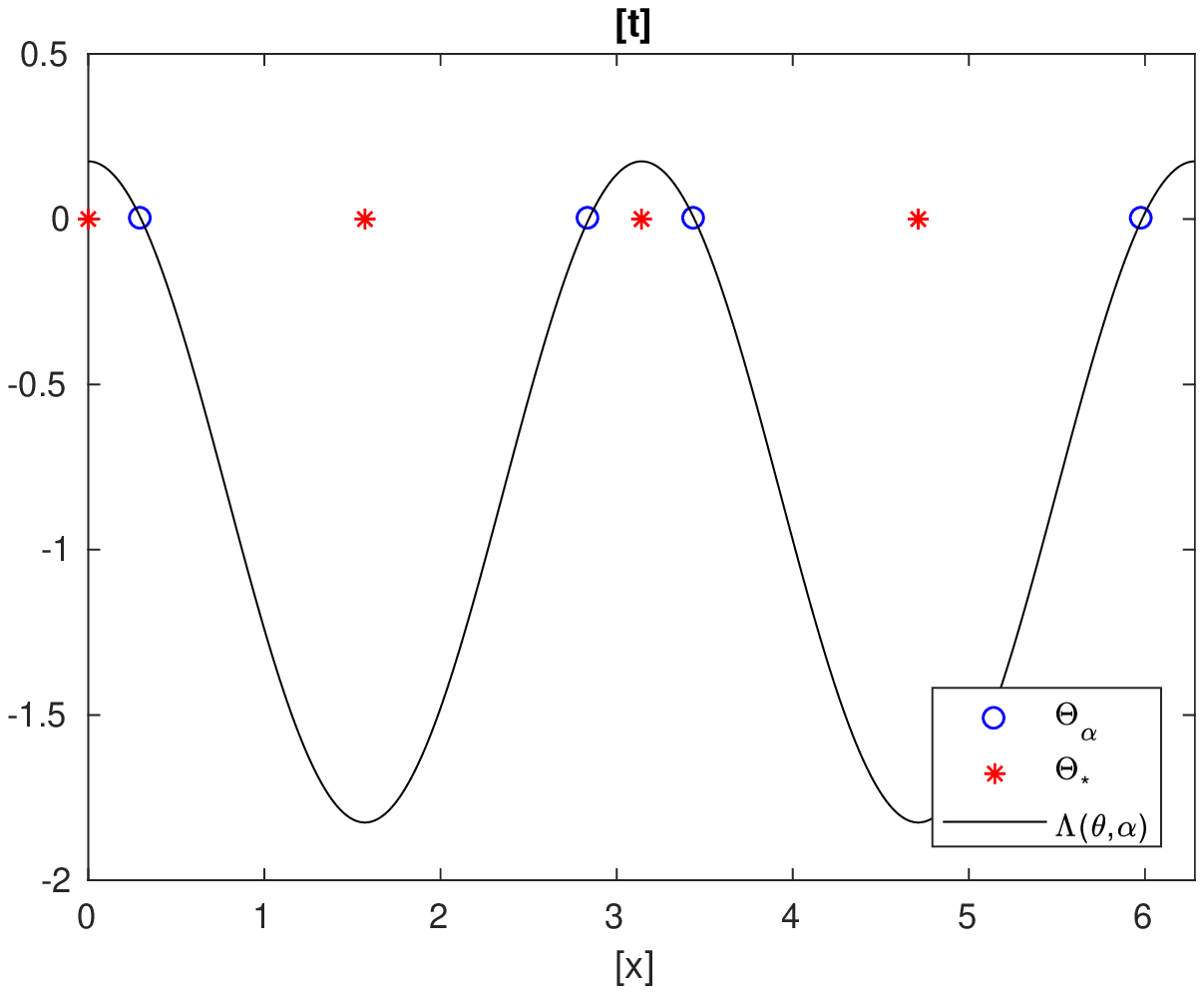}}
		{\psfrag{[x]}{$\theta$}\psfrag{[y]}[Bc][Bc]{}
			\psfrag{[x1]}[Bc][Bc]{$0$}\psfrag{[x2]}[Bc][Bc]{$\pi/2$}\psfrag{[x3]}[Bc][Bc]{$\pi$}\psfrag{[x4]}[Bc][Bc]{$3\pi/2$}\psfrag{[x5]}[Bc][Bc]{$2\pi$}\psfrag{[t]}[Bc][Bc]{$\alpha=0.1,$ $p=4$}\psfrag{[a]}[Bc][Bc]{$\Lambda(\theta,\alpha)$}
			\includegraphics[width=.495\textwidth]{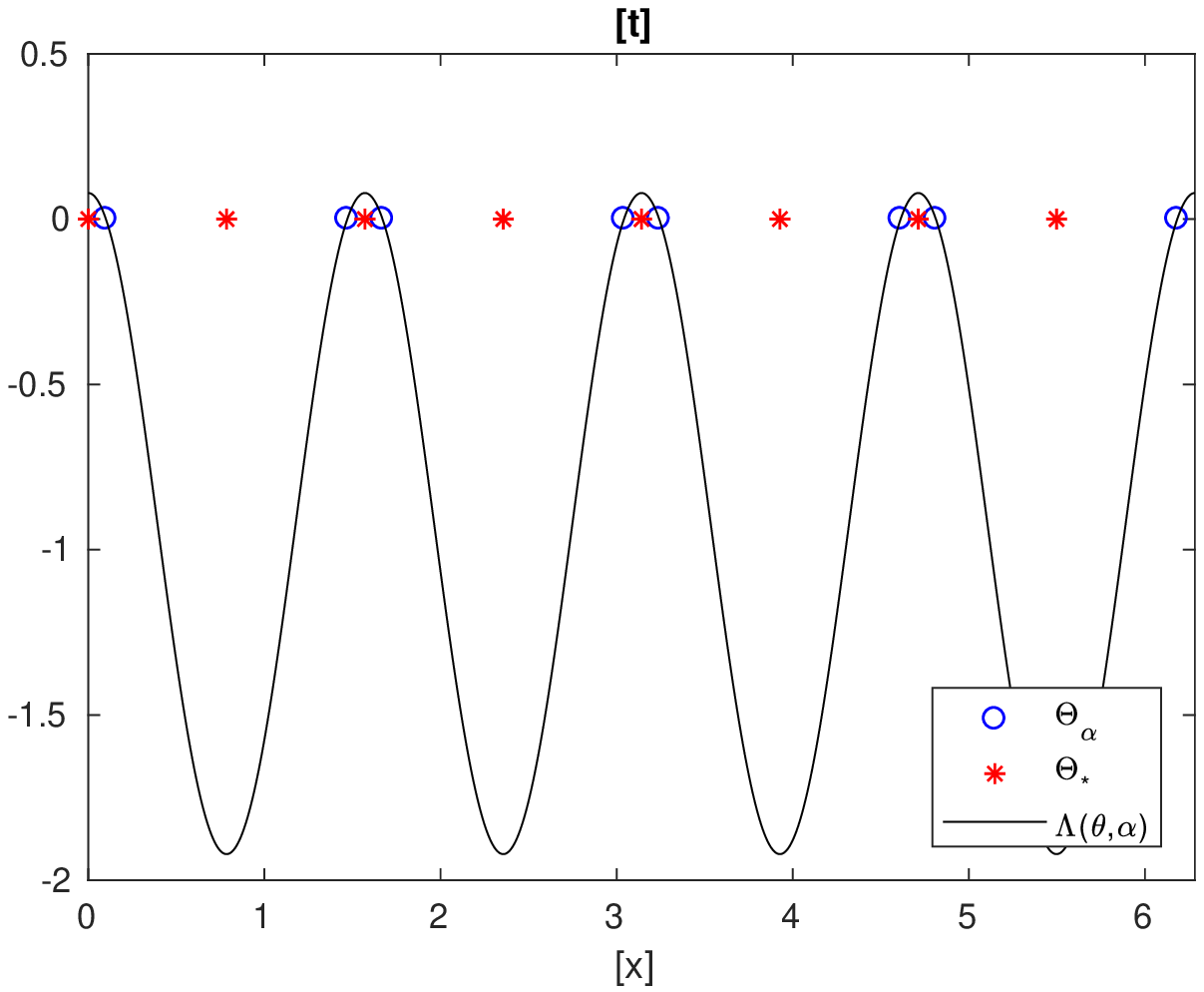}}
		\caption{Two examples of symmetric distribution of $\theta_0\in\Theta_\alpha$ about the points $\theta_*\in\Theta_*$.}\label{fig:symDist}
	\end{figure}
	To see why \eqref{eq:dodgyInequality} holds, we note that the points  $\theta_0\in\Theta_\alpha$ are distributed symmetrically about the points $\theta_*\in\Theta_*$ (as illustrated by Figure \ref{fig:symDist}), and we have specified the condition that $\theta_0$ must be the element of $\Theta_\alpha$ closest to $\theta$, hence it follows that $\theta$ and $\theta_0$ both lie on the same side of $\theta_*$ (in a local sense). For the case $\theta_*=0$ (zero is always an element of $\Theta_*$ due to \eqref{LHoppyLopps2}), it is clear that $|\theta+\theta_0|$ is further from zero than $|\theta+\theta_*|=|\theta|$, whilst still less than $\pi/(2p)$ (given our condition that $\theta_*$ lies close to $\theta$), hence \eqref{eq:dodgyInequality} holds for the case $\theta_*=0$. The case for non-zero $\theta_*\in\Theta_*$ follows similarly, noting that adding linear combinations of elements of $\Theta_*$ to the argument of $|\sin(p(\cdot+\theta_0)/2)|$ does not change its value.
	
	Combining \eqref{Tbound1} with \eqref{eq:dodgyInequality} yields
	\begin{align*}
	|{\Lambda(\theta,\alpha)}|
	&\geq
	{2\left|\sin\left({p}(\theta-\theta_0)/2\right)\right|\cdot\left|\sin\left({p}(\theta_*+\theta_0)/2\right)\right|}\\
	&={2\left|\sin\left({p}(\theta-\theta_0)/{2}\right)\right|\cdot\left|\sin\left({p}(2\theta_*+(\theta_0-\theta_*))/{2}\right)\right|}\\
	&={2\left|\sin\left({p}(\theta-\theta_0)/{2}\right)\right|\cdot\left|\sin(p\theta_*)\cos(p(\theta_0-\theta_*)/2)+\sin(p(\theta_0-\theta_*)/2)\cos(p\theta_*)\right|}.
	\end{align*}
	By the definition \eqref{LHoppyLopps2} of $\Theta_*$, we have that $\sin(p\theta_*)=0$ and $|\cos(p\theta_*)|=1$, hence
	\begin{equation}\label{sinesBeforeBound}
	\left|{\Lambda(\theta,\alpha)}\right|
	\geq{2\left|\sin\left({p}(\theta-\theta_0)/{2}\right)\right|\cdot\left|\sin(p(\theta_0-\theta_*)/2)\right|}.
	\end{equation}
	From the definitions of $\Theta_\alpha$ \eqref{LHoppyLopps} and $\Theta_*$ \eqref{LHoppyLopps2}, it follows that the furthest $\theta$ can be from the nearest $\theta_0$ or $\theta_*$ is $\pi/(2p)$. Hence the argument of both sines of \eqref{sinesBeforeBound} is at most $\pi/4$, so we may use the identity $|\sin(x)|\geq|x/2|$ for $0\leq|x|\leq\pi/4$ (twice) to obtain the lower bound on $|\Lambda(\theta,\alpha)|$ as claimed.
	
	Now we focus on the upper bound. Writing $\theta+\theta_0=(\theta+\theta_*)+(\theta_0-\theta_*)$ and applying elementary trigonometric addition formulae gives
	\[
	\sin\left(\frac{p}{2}(\theta+\theta_0)\right)=
	\sin\left(\frac{p}{2}(\theta+\theta_*)\right)\cos\left(\frac{p}{2}(\theta_0-\theta_*)\right)+\sin\left(\frac{p}{2}(\theta_0-\theta_*)\right)\cos\left(\frac{p}{2}(\theta+\theta_*)\right)
	\]
	which we can bound using the triangle inequality and $\cos(x)\leq1$ to obtain
	\[
	\left|\sin\left(\frac{p}{2}(\theta+\theta_0)\right)\right|\leq
	\left|\sin\left(\frac{p}{2}(\theta+\theta_*)\right)\right|+\left|\sin\left(\frac{p}{2}(\theta_0-\theta_*)\right)\right|.
	\]
	Since $|\sin(px/2)|$ is symmetric about the points $x\in\Theta_*$, it follows that $|\sin\left({p}(\theta+\theta_*)/{2}\right)|=|\sin({p}(\theta-\theta_*)/{2})|$, 
	hence
	\[
	\left|\sin\left(\frac{p}{2}(\theta+\theta_0)\right)\right|\leq
	\left|\sin\left(\frac{p}{2}(\theta-\theta_*)\right)\right|+\left|\sin\left(\frac{p}{2}(\theta_0-\theta_*)\right)\right|.
	\]
	We may use this to bound \eqref{Tbound1}, at which point the identity $|\sin(x)|\leq |x|$ may be used to obtain
	\[
	\left|\Lambda(\theta,\alpha)\right|\leq
	2\frac{p}{2}|\theta-\theta_0|\frac{p}{2}\left(|\theta-\theta_*|+ |\theta_0-\theta_*|\right),
	\]
	finally we use the triangle inequality to write $|\theta-\theta_*|\leq|\theta-\theta_0|+|\theta_0-\theta_*|$ and obtain the upper bound as claimed.\qed
\end{proof}

The upper bound of Lemma~\ref{HatBound} can be used to identify and quantify cases where the error from the solver $\mathcal{P}$ becomes significantly exaggerated by the term $1/\Lambda(\theta,\alpha)$. Conversely, we observe that $|~\theta_0~-~\theta_*~|~\ll~1$ is not a sufficient condition for instability, but when coupled with the additional condition of $|\theta-\theta_0|\ll1$ we observe a stronger instability. These instabilities are demonstrated by the blow-up {of the relative error shown} in Figure~\ref{PS_err}. 
{It is also clear from Figure \ref{PS_err} that increasing the accuracy of the solver $\mathcal{P}$, by increasing the number of degrees of freedom (DOFs) of the underlying method of fundamental solutions, this will reduce the region of instability.}
However, for any degree of accuracy there will always be values of $\theta$ at which the error is arbitrarily large, and as we shall see, increasing DOFs may be an unnecessarily expensive approach to reducing the width of the unstable region.

\newcommand{\relErrorNaive}{${\left|\mathcal{E}_{\mathcal{P}}D(\theta,\alpha)-D(\theta,\alpha)\right|}/{\|D(\cdot,\alpha)\|_{L^2(0,2\pi)}}$}

\begin{figure}[hbpt]
\centering	
{\psfrag{[x]}[bl][bl]{$\theta$}\psfrag{[y]}[Bc][Bc]{\relErrorNaive}
\psfrag{[t]}[Bc][Bc]{$\alpha=1$}\psfrag{[t1]}{$0,$}\psfrag{[t2]}{$\alpha,$}\psfrag{[t3]}[Bc][Bc]{$\pi-\alpha,$}\psfrag{[t4]}[Bc][Bc]{$\pi,$}\psfrag{[t5]}[Bc][Bc]{$\pi+\alpha,$}\psfrag{[t6]}[Bc][Bc]{$2\pi-\alpha,$}\psfrag{[t7]}{$2\pi$}\includegraphics[width=.9\textwidth]{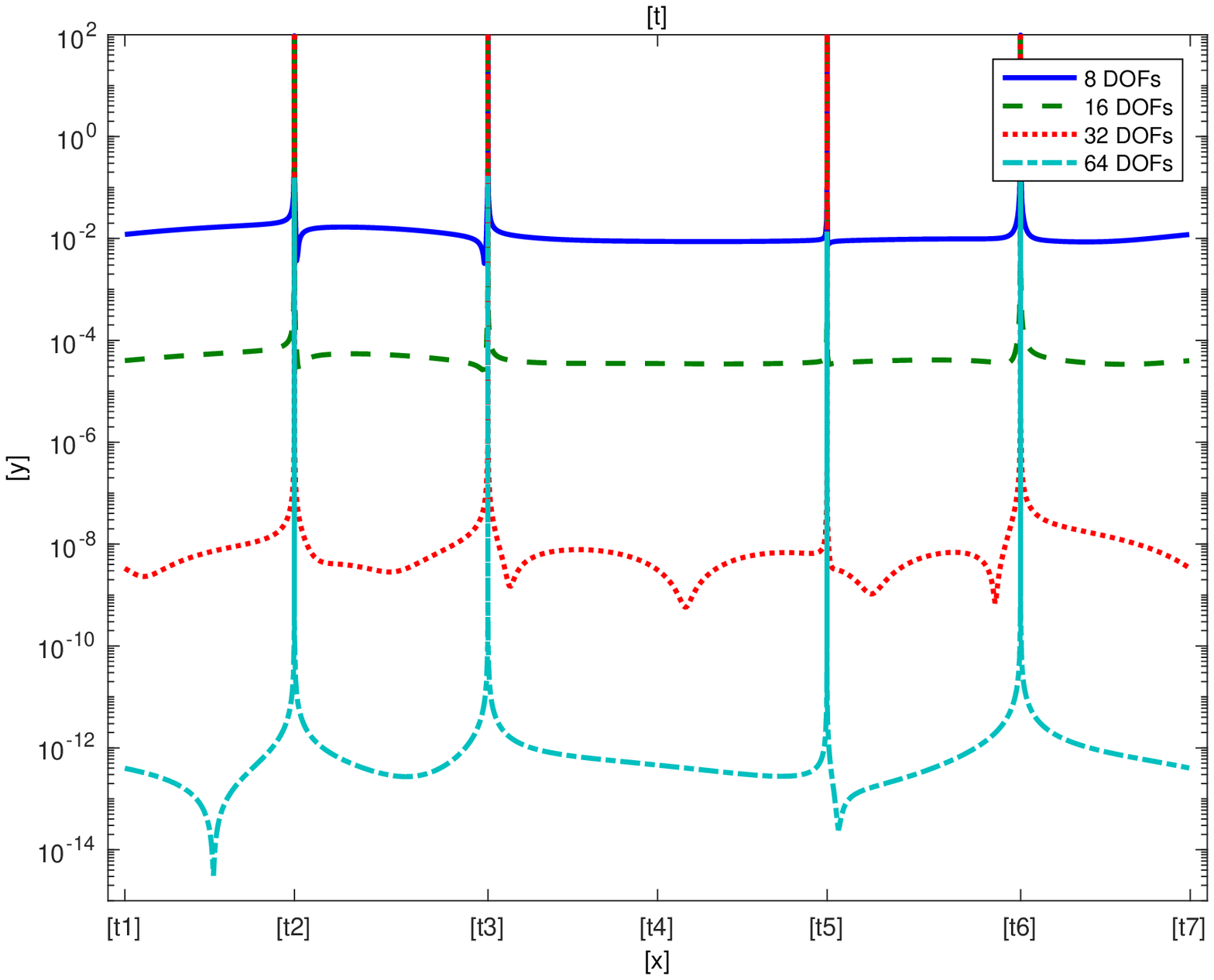}} 
{\psfrag{[x]}[bl][bl]{$\theta$}\psfrag{[y]}[Bc][Bc]{\relErrorNaive }\psfrag{[t]}[Bc][Bc]{$\alpha=0$}\psfrag{[t1]}{$0$}\psfrag{[t2]}{$\pi/2$}\psfrag{[t3]}[Bc][Bc]{$\pi$}\psfrag{[t4]}[Bc][Bc]{$3\pi/4$}\psfrag{[t5]}[Bc][Bc]{$2\pi$}\includegraphics[width=.9\textwidth]{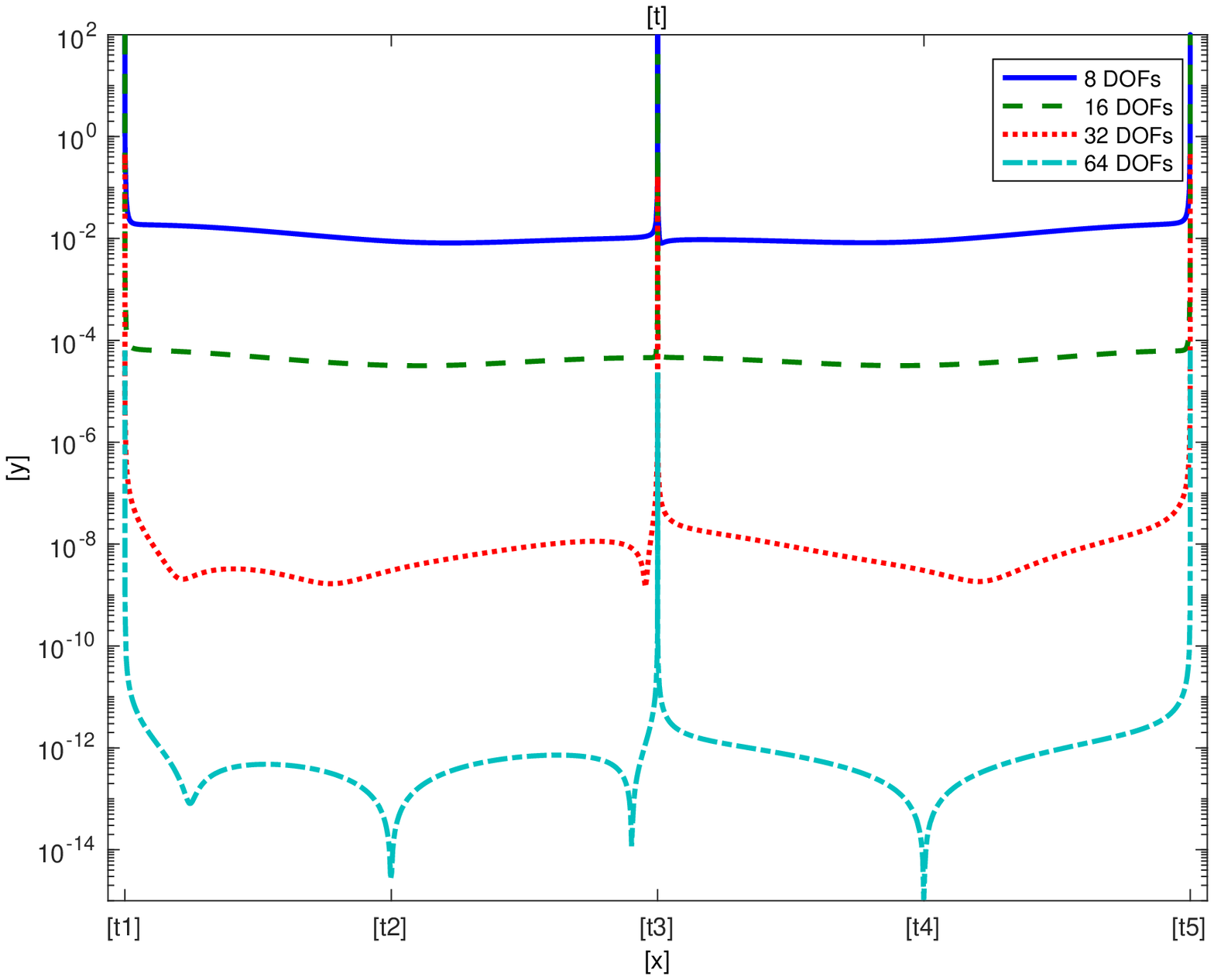}}
\caption{Two plots showing relative error in the naive embedding approximation \eqref{embP0}, for the problem of scattering by a square of side length 1, wavenumber $k=1$, with $\alpha=1$ (not close to $\theta_*\in\Theta_*${, upper panel}), and $\alpha=0\in\Theta_*$ {(lower panel)}. 
Note that there are half as many unstable points for the case $\alpha=0$. 
{Different curves corresponds to different numbers of degrees of freedom in the solution of each of the $M$ canonical problems with MPSpack.}
}
\label{PS_err}
\end{figure}

\section{Numerically stable implementation}\label{stability}
We now 
reformulate the naive implementation \eqref{embP0}, and develop a numerically stable approximation for any $(\theta,\alpha)\in[0,2\pi)^2$. 
{In \S\ref{sec:SVTaylor} we introduce a correction to \eqref{embP0} for $\theta$ close to $\Theta_\alpha$ and far from  $\Theta_0$.
In \S\ref{MV1} we introduce another correction to \eqref{embP0} for $\theta$ close to both $\Theta_\alpha$ and $\Theta_0$.
In \S\ref{s:CEA} all corrections are combined.
}

\subsection{Single variable Taylor expansion}\label{sec:SVTaylor}

{We introduce two small parameters $\mbox{tol}_1\ll1$ and $\mbox{tol}_2\ll1$.}
To address the instability that occurs when
\begin{equation}
|\theta-\theta_0| < \mbox{tol}_1 {\quad \text{and} \quad \mbox{tol}_2 \leq |\theta_*-\theta_0|,}
\label{eqn:theta_tol}
\end{equation}
where $\theta_0$ {and $\theta_*$ are as in Lemma~\ref{HatBound}}, we begin by expanding~(\ref{emb0}) as a Taylor series about $\theta_0$.  {Convergence of this Taylor series is guaranteed by the analyticity of the far-field coefficient (e.g., \cite[Theorem~2.6]{CoKr:13})}. From~\eqref{emb0} it follows that 
\[\sum_{m=1}^\McansChSix B_m(\alpha)\hat{D}(\theta_0,\alpha_m)=\Lambda(\theta_0,\alpha) D(\theta_0,\alpha) = 0, \]
since $\theta_0\in\Theta_\alpha$, hence
the exact representation becomes
\begin{equation}\label{embT}
{D}(\theta,\alpha)=\frac{\theta-\theta_0}{{\Lambda}{(\theta,\alpha)}}{\sum_{m=1}^{\McansChSix}B_m(\alpha)\sum_{n=1}^\infty\frac{(\theta-\theta_0)^{n-1}}{n!}\pdiv{^{n}\hat{D}}{\theta^{n}}(\theta_0,\alpha_m)}. 
\end{equation}
Truncating after $\MTayChSix$ terms, and recalling~(\ref{embP0}), we then define a new approximation
\begin{equation}\label{embPT}
D(\theta,\alpha)\approx\mathcal{E}^{0}_{\mathcal{P}}D(\theta,\alpha; \theta_0,\MTayChSix):=
\frac{\theta-\theta_0}{{\Lambda(\theta,\alpha)}}{\sum_{m=1}^\McansChSix b_m(\alpha)\sum_{n=1}^{\MTayChSix}\frac{(\theta-\theta_0)^{n-1}}{n!}\pdiv{^{n}}{\theta^{n}}\mathcal{P}\hat{D}(\theta,\alpha_m)}.
\end{equation}

Whilst we have introduced an additional source of error in the form of the Taylor remainder, this approximation has advantages over~(\ref{embP0}) when $|\theta-\theta_0|$ is small. 
\newpage
We now have
\begin{align*}
&|D(\theta,\alpha) - \mathcal{E}^{0}_{\mathcal{P}}D(\theta,\alpha; \theta_0,\MTayChSix)| = \\
&\!\left|\frac{\theta-\theta_0}{{\Lambda}{(\theta,\alpha)}} \right|
\left|\sum_{m=1}^{\McansChSix}\sum_{n=1}^{\MTayChSix}\frac{(\theta-\theta_0)^{n-1}}{n!} 
\Big(B_m(\alpha)\pdiv{^{n}\hat{D}}{\theta^{n}}(\theta_0,\alpha_m) - b_m(\alpha)\pdiv{^{n}}{\theta^{n}}\mathcal{P}\hat{D}(\theta,\alpha_m)\Big)\! 
+\mathcal{R}^{0}_{\mathcal{P}}D(\theta,\alpha; \theta_0,\MTayChSix)\right|,
\end{align*}
where
\[ \mathcal{R}^{0}_{\mathcal{P}}D(\theta,\alpha; \theta_0,\MTayChSix) = \sum_{m=1}^{\McansChSix} B_m(\alpha) \sum_{n=\MTayChSix+1}^\infty\frac{(\theta-\theta_0)^{n-1}}{n!}\pdiv{^{n}\hat{D}}{\theta^{n}}(\theta_0,\alpha_m). \]
It follows immediately from Lemma~\ref{HatBound} that 
\[
\left|\frac{\theta-\theta_0}{\Lambda(\theta,\alpha)}\right|\leq\frac{8}{p^2|\theta_*-\theta_0|},\quad\text{for }\theta_*\not\in\Theta_*,
\]
which compares favourably to~\eqref{instab} 
when $\theta$ is close to $\theta_0\in\Theta_\alpha$, provided that $\theta_0$ is sufficiently far from any $\theta_*\in\Theta_*$ {(which is enforced by {the second condition in} \eqref{eqn:theta_tol})}

\begin{figure}[hbpt]
\centering
\subfloat[]
{\psfrag{[x]}{$\theta$}\psfrag{[y]}{$\alpha$}\includegraphics[width=.5\textwidth]{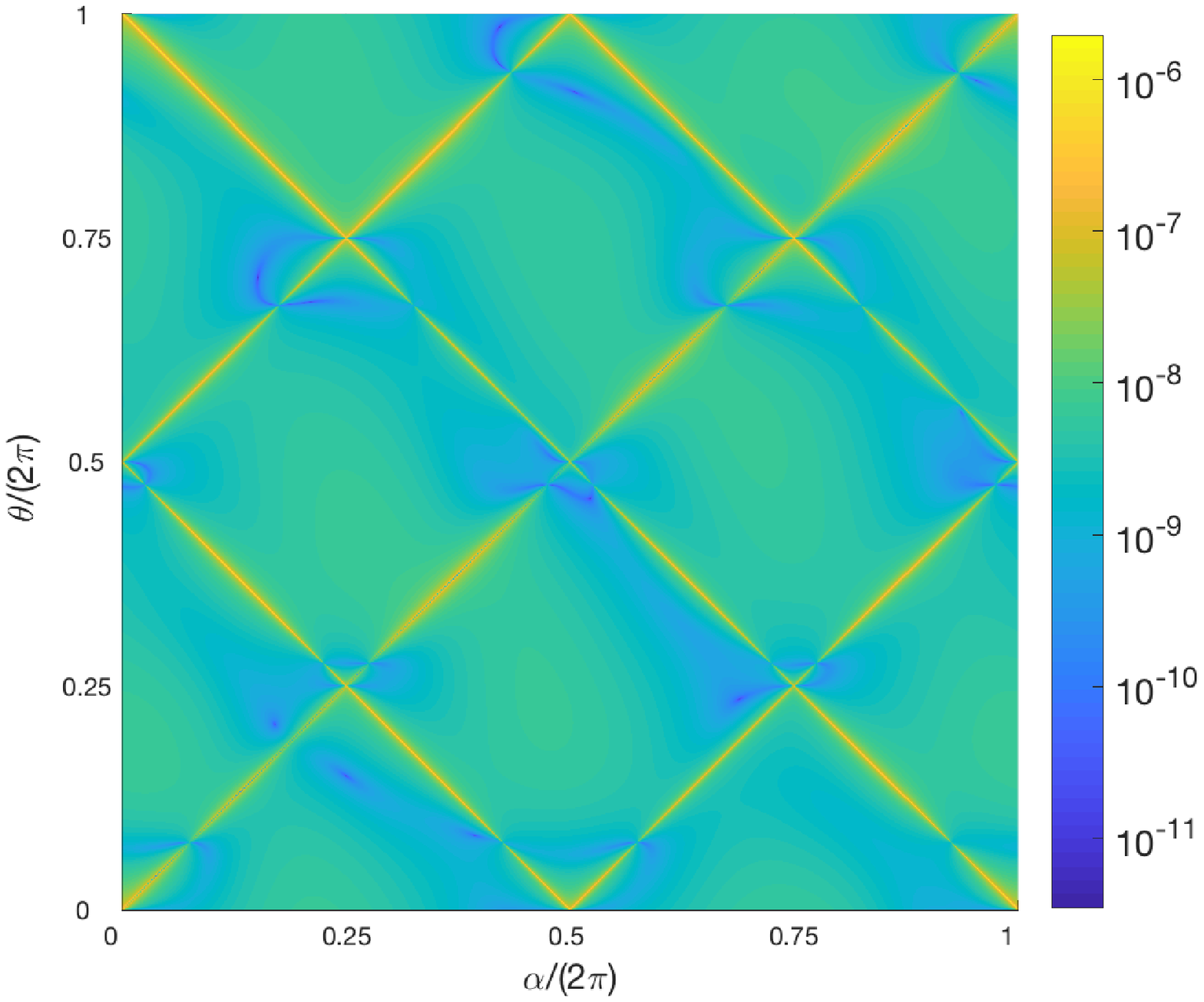}} 
\subfloat[]
{\psfrag{[x]}{$\theta$}\psfrag{[y]}{$\alpha$}\includegraphics[width=.5\textwidth]{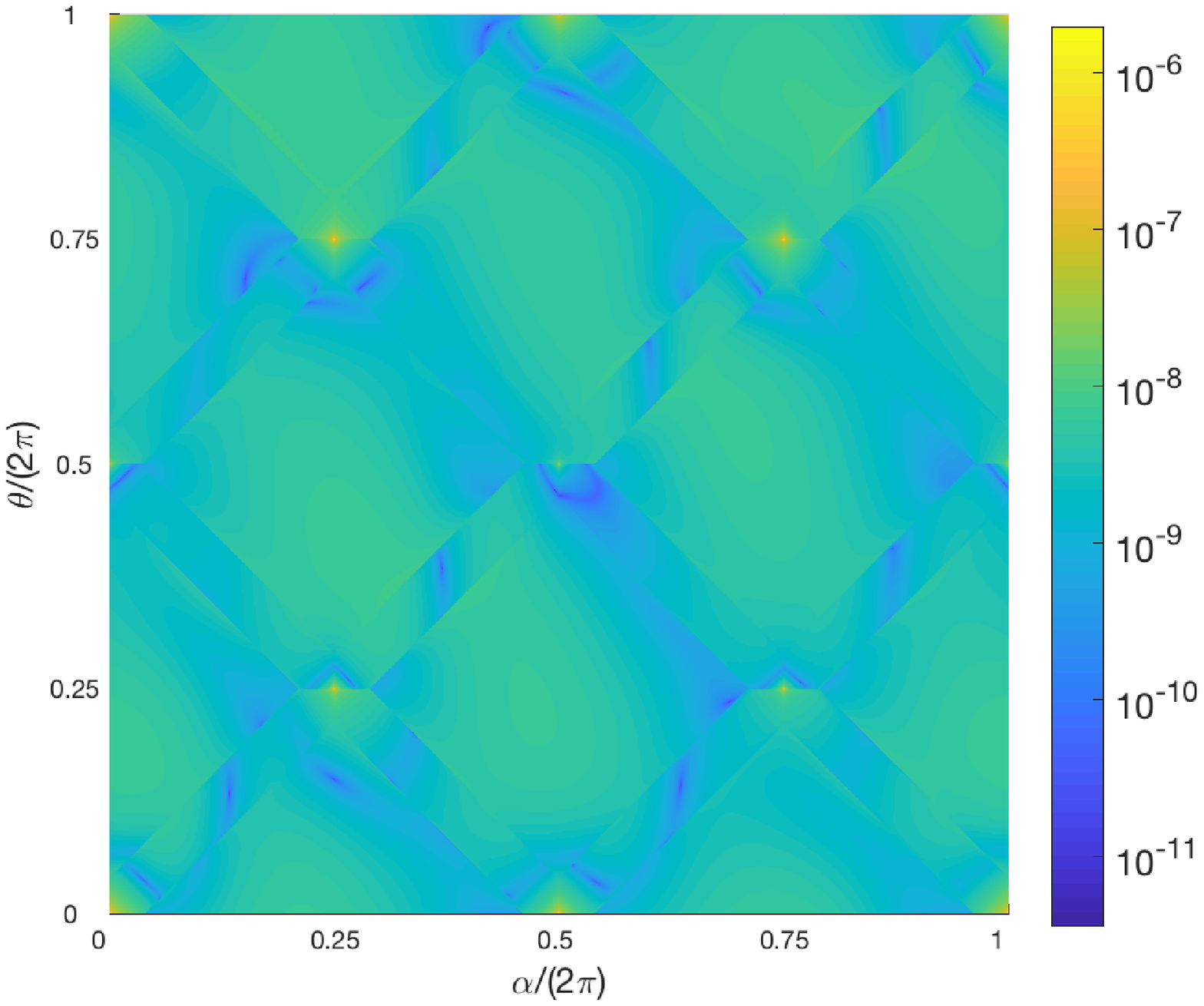}}
\caption{Two plots of relative error for different implementation of embedding formulae for {the problem of scattering by a square of side length 1, wavenumber $k=1$}. Plot (a) depicts accuracy for a naive implementation \eqref{embP0}, whilst Plot (b) depicts accuracy using an $\MTayChSix=10$ degree Taylor expansion approximation \eqref{embPT}, with $\mbox{tol}_1 = 0.25$.
}\label{fig:heat}
\end{figure}

To illustrate the benefits of this approach, in Figure~\ref{fig:heat} we show 
\begin{equation}
\frac{\left|D(\theta,\alpha)-\mathcal{E}_{\mathcal{P}}D(\theta,\alpha)\right|}{\|D\|_{L^2([0,2\pi)^2)}}
\quad\text{and}\quad
\frac{\left|D(\theta,\alpha)-\mathcal{E}_{\mathcal{P}}^0D(\theta,\alpha; \theta_0, 10)\right|}{\|D\|_{L^2([0,2\pi)^2)}}
\end{equation}
in plots (a) and (b), respectively.  Here, we choose $\MTayChSix=10$, and use \eqref{embPT} for all values of $\theta$ such that \eqref{eqn:theta_tol} holds, with $\mbox{tol}_1 = 0.25$.
The reference solution for these and all subsequent results was computed using MPSpack with a large number of degrees of freedom. Whilst the relative error in both cases appears to peak below $10^{-8}$, this is only due to a lack of granularity; these plots contain $1000\times1000$ points but close to the unstable points the relative error is observed to become arbitrarily large for the case $\alpha=0$ in Figure~\ref{PS_err}, which actually focuses on a point which remains unstable in Figure \ref{fig:heat}(b). As Figure~\ref{PS_err} represents the error (against $\theta$) for fixed $\alpha$, it may be considered a vertical cross-section of Figure~\ref{fig:heat}(a). 
We observe that whilst most of the instabilities are handled in Figure \ref{fig:heat}(b) using \eqref{embPT}, there are points at which instabilities still exist. Further care is required to ensure that rounding errors do not cause instability when implementing \eqref{embPT}; a simple fix for this is presented in {\S\ref{sec:denomFix}}. 
However, it is important to note that with a naive implementation of \eqref{embP0}, for every $\alpha\in[0,2\pi)$ there are at least two $\theta_0\in\Theta_\alpha$ for which the approximation is unstable, totalling an infinite number of unstable points $(\theta,\alpha)$ in $[0,2\pi)^2$ (these form the diagonal lines of higher error in Figure \ref{fig:heat}(a)). In contrast when using the approximation \eqref{embPT}, instabilities occur only close to $\theta\in\Theta_*$ for each $\alpha\in\Theta_*$, so \eqref{embPT} reduces the instabilities to a finite set of points, also in $\Theta_*\times\Theta_*$ (these are the points in Figure \ref{fig:heat}(b) at which the lines in Figure \ref{fig:heat}(a) cross).

\subsection{Two-variable Taylor expansion}\label{MV1}
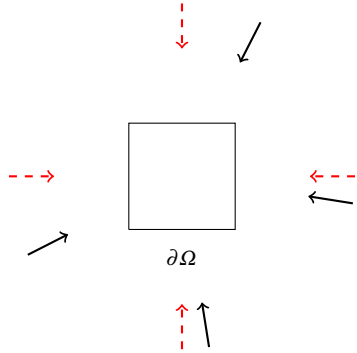
\begin{figure}[hbpt]
	\centering
	\begin{tikzpicture}\draw (-0.707107,-0.707107) -- (0.707107,-0.707107) -- (0.707107,0.707107) -- (-0.707107,0.707107) -- (-0.707107,-0.707107);
	\draw[red,thick,dashed,->] (-2.300000,0.000000) -- (-1.700000,0.000000);
	\draw[red,thick,dashed,->] (-0.000000,-2.300000) -- (-0.000000,-1.700000);
	\draw[red,thick,dashed,->] (2.300000,-0.000000) -- (1.700000,-0.000000);
	\draw[red,thick,dashed,->] (0.000000,2.300000) -- (0.000000,1.700000);
	\draw[thick,->] (-2.049315,-1.044178) -- (-1.514711,-0.771784);
	\draw[thick,->] (0.359799,-2.271683) -- (0.265939,-1.679070);\draw[thick,->] (2.271683,-0.359799) -- (1.679070,-0.265939);\draw[thick,->] (1.044178,2.049315) -- (0.771784,1.514711);\draw (-0.000000,-1.082107) node{$\partial\Omega$};\end{tikzpicture}
	\caption{Diagram of incident angles chosen for a square, when implementing \eqref{shiTaylor}. Incident angles $\alpha_m\in A_\McansChSix\cap\Theta_*$ are {\color{red}dashed}, whilst other incident angles $\alpha_m\in A_\McansChSix\setminus\Theta_*$ are chosen randomly, in accordance with the conditions suggested in Remark \ref{Mchoice2}. }
	\label{MV_A_diag}
\end{figure}

{We now consider the case where
\begin{equation}
|\theta-\theta_0| < \mbox{tol}_1 \quad \text{and} \quad  |\theta_*-\theta_0|<\mbox{tol}_2.
\end{equation}
}Until now, we have been able to choose the canonical incident angles $A_\McansChSix$ relatively freely (as discussed in Remark \ref{Mchoice2}). Contrary to the approach for standard implementation of embedding formulae, it can be advantageous to choose these such that $A_\McansChSix\subset { \Theta_*}$.
Supposing we are close to a source of instability $(\theta_*,\alpha_*){\in \Theta_*\times\Theta_*}$, which as suggested by Figure \ref{fig:heat}(b) is not regulated by \eqref{embPT}, and that the canonical incident waves have been chosen such that $\theta_*=:\alpha_{m_1}\in A_\McansChSix$ and $\alpha_*=:\alpha_{m_2}\in A_\McansChSix$ (see Figure \ref{MV_A_diag}). We can then use the reciprocity relation (Theorem \ref{TH:recip}) to obtain a first order multi-variate Taylor series approximation in terms of quantities that require only the canonical solutions,\newpage
\
\begin{align}
D(\theta,\alpha)&\approx
D(\theta_*,\alpha_*)+(\theta-\theta_*)\pdiv{D}{\theta}(\theta_*,\alpha_*)+
(\alpha-\alpha_*)\pdiv{D}{\alpha}(\theta_*,\alpha_*)
\nonumber\\
&=
D(\theta_*,\alpha_{m_2})+(\theta-\theta_*)\pdiv{D}{\theta}(\theta_*,\alpha_{m_2})+
(\alpha-\alpha_*)\pdiv{D}{\theta}(\alpha_*,\alpha_{m_1})
\nonumber
\\
&\approx\mathcal{P}D(\theta_*,\alpha_{m_2})+(\theta-\theta_{*})\pdiv{\mathcal{P}D}{\theta}(\theta_*,\alpha_{m_2})+
(\alpha-\alpha_*)\pdiv{\mathcal{P}D}{\theta}(\alpha_*,\alpha_{m_1})\nonumber
\\
&=:\mathcal{E}_{\mathcal{P}}^*D(\theta,\alpha; \theta_*,\alpha_*),
\label{shiTaylor}
\end{align}
where the derivative of each far-field coefficient {can be computed} from the representation \eqref{FFrep}. 
We shall incorporate \eqref{shiTaylor} as one of several approximations in a method in \S\ref{s:CEA}.

\begin{remark}\label{MVh}
	For simplicity of notation, we just show in \eqref{shiTaylor} the first order two-variable Taylor expansion.  In principle this can be extended to higher orders, but this introduces some challenges in computing the cross derivative terms.  See \cite[Appendix~C.3]{Gi:17} for further details.
\end{remark}

\subsection{Combined approach}\label{s:CEA}
We are now able to construct a numerically stable (over all $[0,2\pi)^2$) 
method using an appropriate combination of the approximations \eqref{embP0} (with L'Hopital's rule when $\theta\in\Theta_{\alpha}$), \eqref{embPT} and \eqref{shiTaylor}. 
\begin{definition}[Combined embedding approximation]\label{def:CEA}
\; Given thresholds $0<\mbox{tol}_1\ll1$ and $0<\mbox{tol}_2\ll1$, we choose the canonical set of incident waves $A_\McansChSix$ such that $\Theta_*\subset A_\McansChSix$. We define the combined embedding approximation by
\begin{align}\label{CEEAdef}
&\mathcal{E}_{\mathcal{P}}^{\circledast}D(\theta,\alpha;\MTayChSix):=\\
&\left\{\begin{aligned}
&\left[\sum_{m=1}^\McansChSix b_m(\alpha)\frac{\partial\mathcal{P}\hat{D}}{\partial \theta}(\theta,\alpha_m)\right]\Big/&&\hspace{-4mm}\left[\pdiv{\Lambda}{\theta}(\theta,\alpha)\right],\qquad \text{if }\theta\in\Theta_{\alpha}\text{ and }\min_{\theta_*\in\Theta_*}|\theta-\theta_*|>\mbox{tol}_2, \text{ else}\\
&\left[\sum_{m=1}^\McansChSix b_m(\alpha)\frac{\partial^2\mathcal{P}\hat{D}}{\partial \theta^2}(\theta,\alpha_m)\right]\Big/&&\hspace{-4mm}\left[\pdiv{^2\Lambda}{\theta^2}(\theta,\alpha)\right], \qquad\text{if }\theta\in\Theta_{\alpha}\text{ and }\theta\in\Theta_*, \text{ else}\\
&\mathcal{E}_{\mathcal{P}}^{0}D(\theta,\alpha; \theta_0,\MTayChSix),\quad 
&&\text{if } \theta_0\in\Theta_\alpha:0<|\theta-\theta_0|\leq \mbox{tol}_1\text{ and }\min_{\theta_*\in\Theta_*}|\theta-\theta_*|>\mbox{tol}_2, \text{ else}\\
&\mathcal{E}_{\mathcal{P}}^{*}D(\theta,\alpha; \theta_*,\alpha_*),
&&\text{if }|\theta-\theta_0|<\mbox{tol}_2\text{ and }|\theta_0-\theta_*|<\mbox{tol}_2,
\text{ for some }\theta_*\text{ and }\alpha_*\in\Theta_{*}, \\
&\mathcal{E}_{\mathcal{P}}D(\theta,\alpha),&&\text{otherwise},
\end{aligned}\right.
\end{align}
for $(\theta,\alpha)\in[0,2\pi)^2$, where $\theta_0$ is chosen to be the element of $\Theta_\alpha$ closest to $\theta$, and $\theta_*$ is chosen to be the element of $\Theta_*$ closest to $\theta_0$ with $\alpha_*\in\Theta_*$ such that $\theta_*\in\Theta_{\alpha_*}$.
{(Recall that $\mathcal{E}_{\mathcal{P}}^{0}$, $\mathcal{E}_{\mathcal{P}}^{*}$ and $\mathcal{E}_{\mathcal{P}}$ were defined in \eqref{embPT}, \eqref{shiTaylor} and \eqref{embP0}, respectively.)}
\end{definition}
Results for the combined approach are depicted in Figure~\ref{fixed_fig}. The relative error is determined by comparing against the MPSpack approximation to the far-field, $\mathcal{P}D$. The error observed is greater for $\alpha=0$ than for Figure~\ref{fig:heat}(a), since we have used the first order expansion \eqref{shiTaylor} instead of the tenth order expansion used for Figure~\ref{fig:heat}(a). When compared against the results of Figure~\ref{PS_err}, a significant improvement is visible.  In the lower plot in Figure~\ref{fixed_fig}, we show a close up for $\alpha=0$ around $\theta=\pi$, demonstrating that the unbounded error observed in Figure~\ref{PS_err} is no longer present.

\newcommand{\relErrorCombined}{${\left|\mathcal{E}^\circledast_{\mathcal{P}}D(\theta,\alpha)-D(\theta,\alpha)\right|}/{\|D(\cdot,\alpha)\|_{L^2(0,2\pi)}}$}

\begin{figure}[hbpt]
	\footnotesize
	\centering
	{\psfrag{[x]}[bl][bl]{$\theta$}
		\psfrag{[y]}[Bc][Bc]{\relErrorCombined}
		\psfrag{[t]}{$\alpha=1$}\psfrag{[t1]}{$0,$}\psfrag{[t2]}{$\alpha,$}\psfrag{[t3]}[Bc][Bc]{$\pi-\alpha,$}\psfrag{[t4]}[Bc][Bc]{$\pi,$}\psfrag{[t5]}[Bc][Bc]{$\pi+\alpha,$}\psfrag{[t6]}[Bc][Bc]{$2\pi-\alpha,$}\psfrag{[t7]}{$2\pi$}
		\includegraphics[width=.6\textwidth]{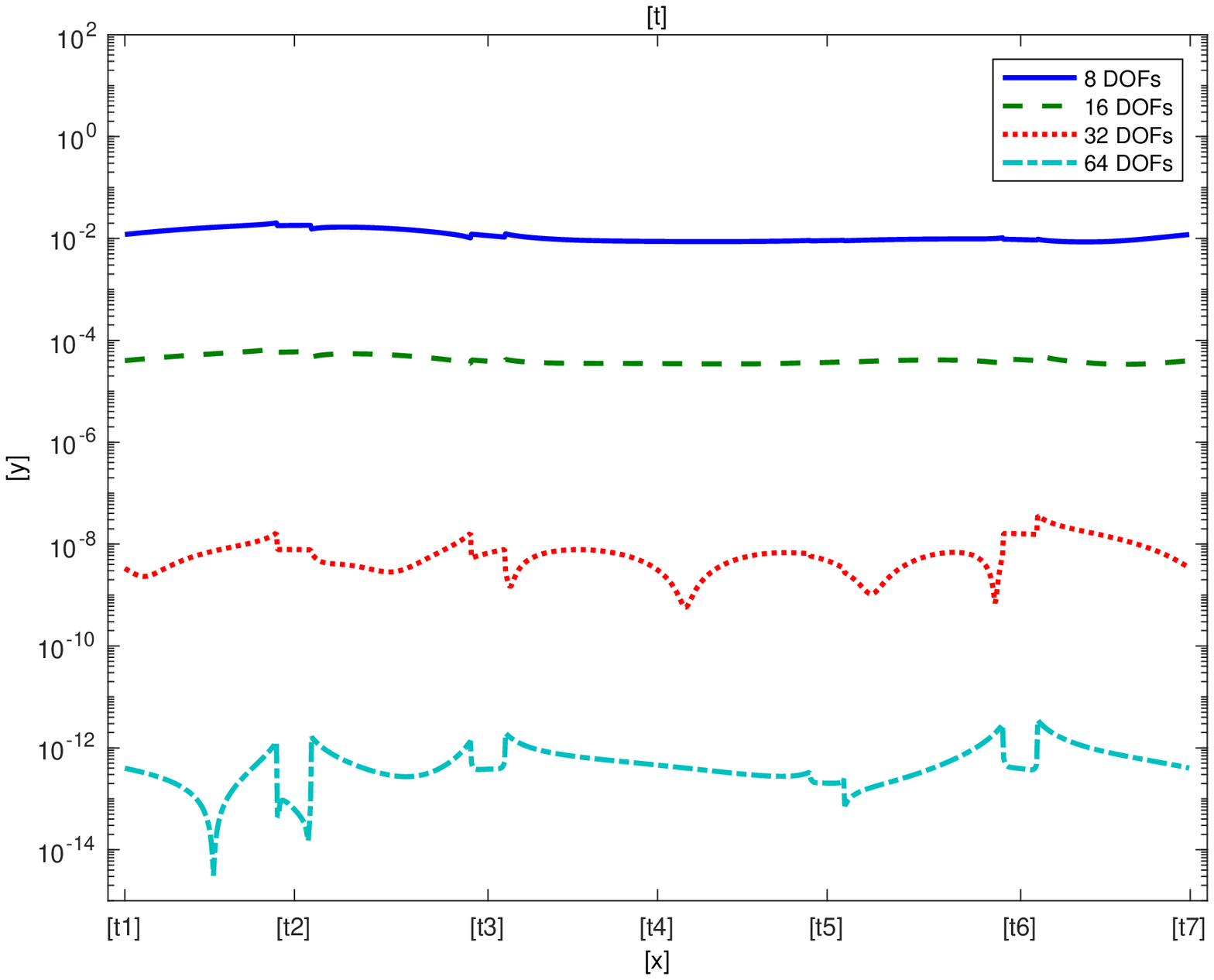}} 
	\centering
	{\psfrag{[x]}[bl][bl]{$\theta$}
		\psfrag{[y]}[Bc][Bc]{\relErrorCombined}
		\psfrag{[t]}{$\alpha=0$}\psfrag{[t1]}{$0$}\psfrag{[t2]}{$\pi/2$}\psfrag{[t3]}[Bc][Bc]{$\pi$}\psfrag{[t4]}[Bc][Bc]{$3\pi/4$}\psfrag{[t5]}[Bc][Bc]{$2\pi$}\includegraphics[width=.6\textwidth]{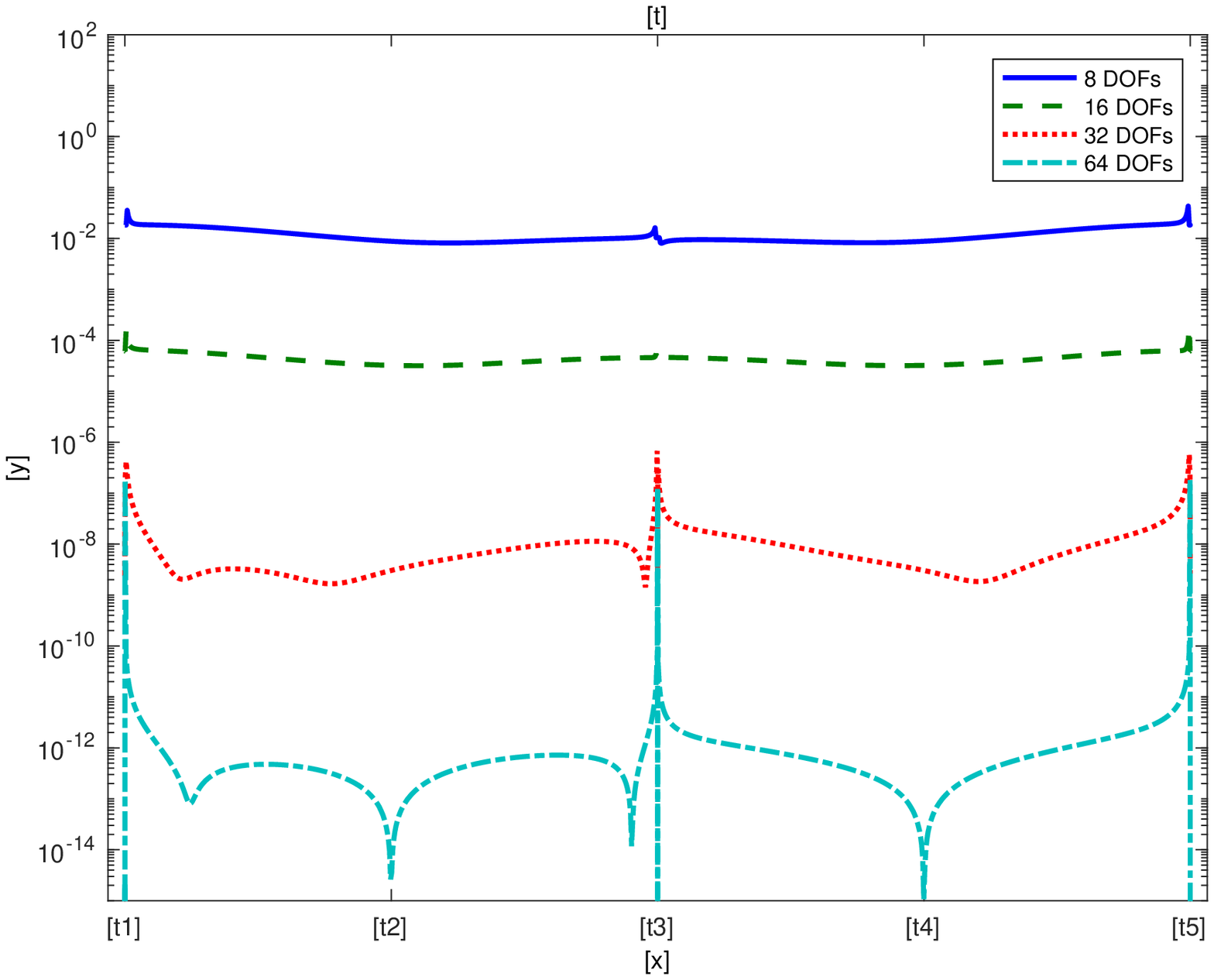}}

	\centering
	{\psfrag{[x]}[bl][bl]{$\theta$}\psfrag{[x1]}[bl][bl]{$0.999\pi$} \psfrag{[x2]}[bl][bl]{$\pi$} \psfrag{[x3]}[bl][bl]{$1.001\pi$}
		\psfrag{[y]}[Bc][Bc]{\relErrorCombined}
		
		\psfrag{[t]}[Bc][Bc]{$\alpha=0$, focusing on $\pi\in\Theta_*$}\psfrag{[t1]}{$0$}\psfrag{[t2]}{$\pi/2$}\psfrag{[t3]}[Bc][Bc]{$\pi$}\psfrag{[t4]}[Bc][Bc]{$3\pi/4$}\psfrag{[t5]}[Bc][Bc]{$2\pi$}\includegraphics[width=.6\textwidth]{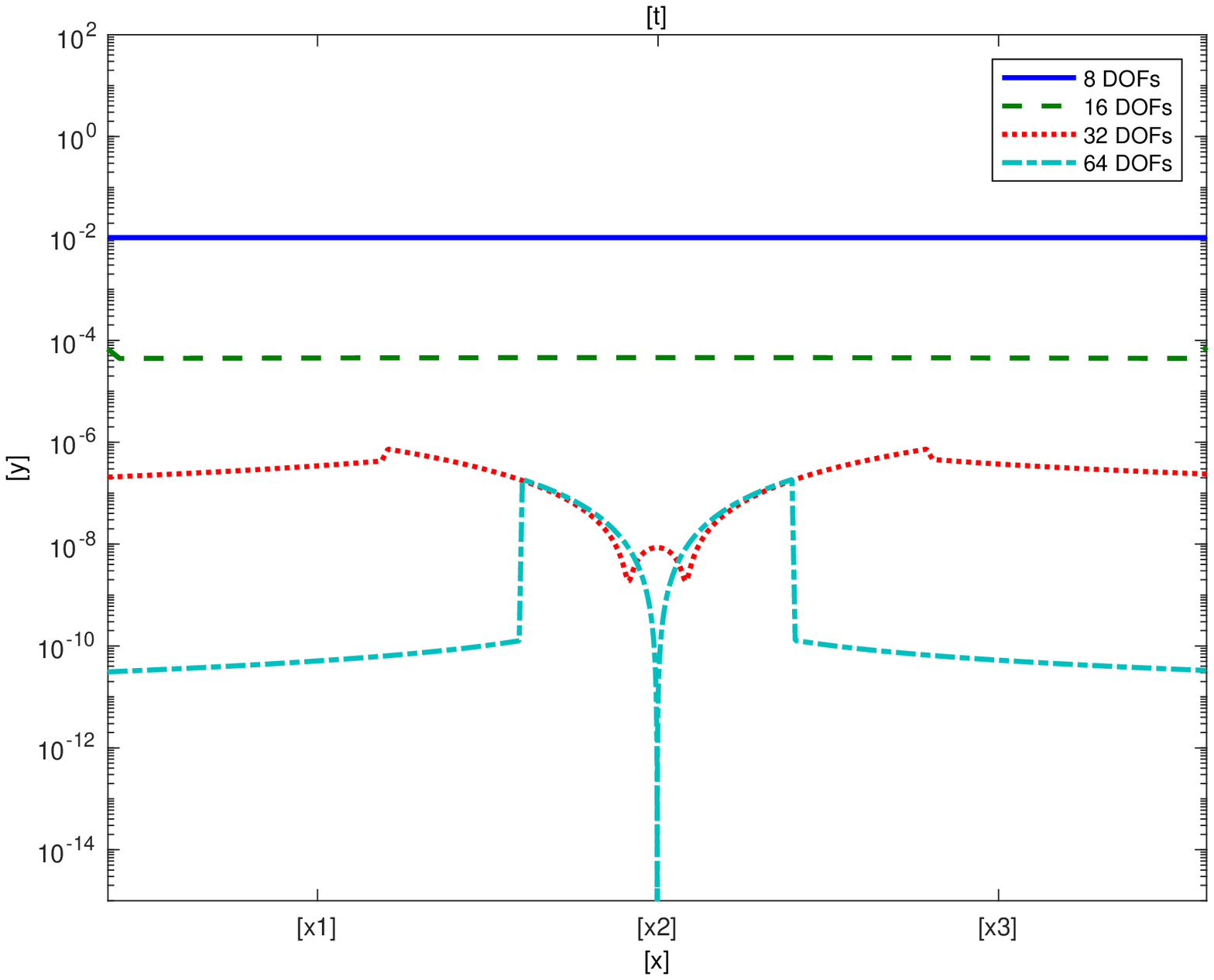}}
	\caption{Relative error for the combined embedding approximation, with $\alpha=1$ (not close to $\theta_*\in\Theta_*$), and $\alpha=0\in\Theta_*$. The lower plot is a close-up for $\alpha=0$ around $\theta=\pi$.}\label{fixed_fig}
\end{figure}

	\section{Practical considerations}\label{sec:imp}
	In this section we briefly introduce two ideas which should be considered for an implementation of the embedding formulae.
	\subsection{Avoiding numerical instabilities in computation of Taylor expansion}\label{sec:denomFix}
	Practical implementation of \eqref{embPT} also requires care when $0<|\theta-\theta_0|\ll 1$, although this is easier to remedy than for the naive implementation \eqref{embP0}. This is because numerical errors in the denominator of
	\[
	\frac{\theta-\theta_0}{\Lambda(\theta,\alpha)}=\frac{\theta-\theta_0}{\cos(p\theta)-(-1)^p\cos(p\alpha)}
	\]
	can cause the value of the total fraction to become significantly inaccurate. Using the identity \eqref{alpha_replace_theta0} and Taylor expanding about $\theta_0$ yields the stable (in a region around $\theta_0$) representation, for $\alpha\not\in\Theta_*$,
	\[
	\frac{\theta-\theta_0}{\Lambda(\theta,\alpha)}=\left({\sum_{n=1}^\infty\frac{(\theta-\theta_0)^{n-1}}{n!}p^n\cos(p\theta_0+n\pi/2)}\right)^{-1}.
	\]
	In practice the sum should be truncated after an appropriate number of terms.
	
	\subsection{Fast algorithm for computing far-field derivatives}\label{sec:fastFF}
	Implementation of \eqref{embPT} also requires computation of derivatives of the far-field coefficient, which can be represented using \eqref{FFrep} as
	\begin{equation}\label{eq:ffDir}
	\pdiv{^{n}D}{\theta^{n}}(\theta,\alpha)=-\int_\Gamma \pdiv{^{n}K}{\theta^{n}}(\theta,\bfy)\pdiv{u_\alpha}{\bfn}(\bfy)\dd{s}(\bfy),
	\end{equation}
	where $K(\theta,\bfy):=\e^{-\imag k [y_1\cos\theta+y_2\sin\theta]}$. We require an efficient method to compute derivatives of the kernel $K$. We can write these derivatives as
	$
	\pdivl{^{n}K}{\theta^{n}}=g_n\cdot K
	$
	where $g_n$ can be defined iteratively,
	\begin{equation}\label{g_1_def}
	\quad
	g_1(\theta,\bfy):=-\imag k[-y_1\sin(\theta)+ y_2\cos(\theta)]
	\end{equation}
	and (by repeated application of the product and chain rules)
	$$
	g_n(\theta,\bfy) = g_{n-1}(\theta,\bfy){g_1}(\theta,\bfy)+\pdiv{g_{n-1}}{\theta}(\theta,\bfy)
	=  \Big({g_1}(\theta,\bfy)+\pdiv{}{\theta}\Big)^ng_1(\theta,\bfy),\qquad\text{for }n\geq2.
	$$
	In practice, the operation which maps $g_{n-1}$ to $g_n$ can be computed
	{as a (banded, rectangular) matrix-vector product for functions that can be expressed as finite linear combinations of complex exponentials $\{\e^{\imag\ell\theta}\}_{\ell\in\Z}$.}
	For example, we can rewrite \eqref{g_1_def} as
	\[
	g_1(\theta,\bfy):=c_{-1}(\bfy)\e^{-\imag\theta}+c_1(\bfy)\e^{\imag\theta},
	\]
	where $c_{-1}(\bfy)=-\imag k(y_2-\imag y_1)/2$ and $c_{1}(\bfy)=-\imag k(y_2+\imag y_1)/2$.  More generally, if a function can be represented exactly using $2n+1$ coefficients $b_{-n},\ldots,b_n$, then 
	the operator $
	{g_1}(\theta,\bfy)+\pdivl{}{\theta}
	$ 
	acting on a function of the form $
	\varphi(\theta; \bfy)=\sum_{\ell=-n}^{n}b_\ell(\bfy)\e^{\imag\ell\theta}
	$ 
	can be represented by (dropping the argument $(\bfy)$ from $c_{-1}$ and $c_{1}$ for notational convenience)
	\begin{equation}\label{eq:GnRep}
	\tens{{G}_n}:=\left[\begin{array}{ccccc}
	c_{-1} & 0 & \cdots &&\\
	-n 		 & c_{-1} & \ddots &&\\
	c_{1}	 & -n +1	& \ddots &&\vdots\\
	0	& \ddots		& \ddots & \ddots &0\\
	\vdots& \ddots & & n-1 &	c_{-1}\\
	& & \ddots& \ddots& n\\
	& & & 0 & c_1\\
	\end{array}\right],
	\end{equation}
	where the $-n,\ldots,n$ terms come from the differentiation and the other terms from the operation of multiplying by $g_1$ (see \cite[Appendix~C.1]{Gi:17} for details).  The fast computation of the far-field kernel follows by repeated multiplication of these matrices
	\begin{equation}\label{eq:gnRep}
	g_n(\theta; \bfy)=\left({g_1}(\theta,\bfy)+\pdiv{}{\theta}\right)^ng_1(\theta; \bfy)
	=
	\Big[\tens{G_n}\times\ldots\times \tens{G_1}[c_{-1}, c_1]^T\Big]^T
	\left[\e^{\imag\ell\theta}\right]_{\ell=-n}^n.
	\end{equation}
	Given $c_{-1}$, $c_1$ and $n$, arrays of $\tens{G_n}$ may be computed very easily (using \eqref{eq:GnRep}), allowing for fast (non-symbolic) computation of far-field derivatives, using~\eqref{eq:gnRep} and subsequently~\eqref{eq:ffDir}.

\begin{remark}[Further error analysis] It follows by the Definition \ref{def:CEA} that the combined embedding approximation avoids the issue of the standard approach, highlighted by Figure~\ref{fig:heat}(a), in which the errors inherited from the solver $\mathcal{P}$ are multiplied by an arbitrarily large number. Some necessary components for a full error analysis, including a fully-explicit upper bound on the truncation error of the Taylor expansions of Definition \ref{s:CEA}, are given in  \cite[Theorem~6.8,~Lemma~6.9]{Gi:17}. However, a full error analysis requires an estimate of each $|B_m(\alpha)-b_m(\alpha)|$ for $m=1,\ldots,M$, which depends on
{the norm of the inverse of the matrix of the system in Assumption \ref{b_m_unique}(i).}
Numerical experiments suggest that this {norm} grows as $k\rightarrow0$, however we do not fully understand how {it} depends on $k$ or on the angles of the canonical incident waves, nor do we understand sufficient conditions for Assumption \ref{b_m_unique} to hold.
\end{remark}

\section{More general incident waves}\label{sec:HergEmbed}
{}

{We now demonstrate how embedding formulae may be used to approximate the far-field coefficient of a far broader class of incident waves than just plane waves, by means of a general formula and some numerical examples. 
In particular, we will be interested in incident waves of the following form (see e.g. \cite[Definition~3.18]{CoKr:13}), which can be thought as continuous linear combinations of plane waves.

\begin{definition}[Herglotz wave functions]\label{hergs}
Given $\gHerg\in L^2(0,2\pi)$, the function 
\[
\uiHerg(\bfx;\gHerg)=\int_0^{2\pi}\gHerg(\alpha)\e^{\imag k \bfx\cdot\bfd_\alpha}\dd{\alpha},\quad\text{for }\bfx\in\R^2,
\]
where $\bfd_\alpha:=-(\cos\alpha,\sin\alpha)$,
is called a \emph{Herglotz wave function} or \emph{Herglotz incident field} with Herglotz kernel $\gHerg\in L^2(0,2\pi)$.
\end{definition}

A second concept which we will find useful is the \emph{far-field map} $\mathcal F_\infty$.
We first extend the scattering boundary value problem \eqref{Helm}--\eqref{SRC} to general incident fields:
if $u^i\in C^\infty(\R^2)$ is an entire solution of the Helmholtz equation, we denote by $u=u^i+u^s$ the solution of the Helmholtz equation \eqref{Helm} in the complement of the scatterer, with $u=0$ on $\partial\Omega$ and such that $u^s$ satisfies the Sommerfeld radiation condition.
Then we write $\mathcal{F}_\infty u^i\in C^\infty(0,2\pi)$ for the far-field coefficient of the scattered field $u^s$.
For example, in terms of familiar plane wave notation \eqref{def:PWinc} and \eqref{FFrep}, we have 
$\mathcal{F}_\infty\uiPW=D(\cdot,\alpha)$.

By Definition \ref{hergs}, the far-field coefficient of a Herglotz {wave} function impinging on $\Omega$ can be computed by integrating the plane wave far-field coefficient against the Herglotz kernel, so it can be approximated using the combined embedding approximation of Definition \ref{CEEAdef}:
\begin{equation}\label{HhergDef}
\big[\mathcal F_\infty \uiHerg(\cdot;\gHerg)\big](\theta)
=\int_0^{2\pi}\gHerg(\alpha)D(\theta,\alpha)\dd{\alpha}
\approx\int_0^{2\pi}\gHerg(\alpha)\mathcal{E}_{\mathcal{P}}^\circledast D\left(\theta,\alpha; \MTayChSix\right)\dd{\alpha}.
\end{equation}
In practice, the integral may be approximated by a quadrature rule, in which $[{w}_i]_{i=1}^{\Nquad}$ and $[\tilde\alpha_i]_{i=1}^{\Nquad}$ are ${\Nquad}\in\mathbb{N}$ suitably chosen weights and nodes sufficient to resolve oscillations and singularities of the integrand to any desired accuracy---the requirement to evaluate the integrand for many incident directions makes this problem particularly well suited to the use of embedding formulae. We may now generalise our combined embedding approximation to any Herglotz wave functions via
\begin{equation}\label{HapproxDef}
\big[\mathcal{E}_{\mathcal{P}}^\circledast 
\mathcal F_\infty \uiHerg(\cdot;\gHerg)\big]\left(\theta,\MTayChSix, ({w}_i)_{i=1}^{\Nquad},(\tilde\alpha_i)_{i=1}^{\Nquad}\right)
:=\sum_{i=1}^{\Nquad}w_i\gHerg(\tilde\alpha_i)\mathcal{E}_{\mathcal{P}}^\circledast D(\theta,\tilde\alpha_i)
\approx
\big[\mathcal F_\infty \uiHerg(\cdot;\gHerg)\big](\theta).
\end{equation}}

\subsection{Numerical example: regular wavefunctions}\label{sec:radwaveEMB}
\begin{figure}
\centering\psfrag{[x]}{$x_1$}\psfrag{[y]}{$x_2$}
\includegraphics[width=0.8\textwidth]{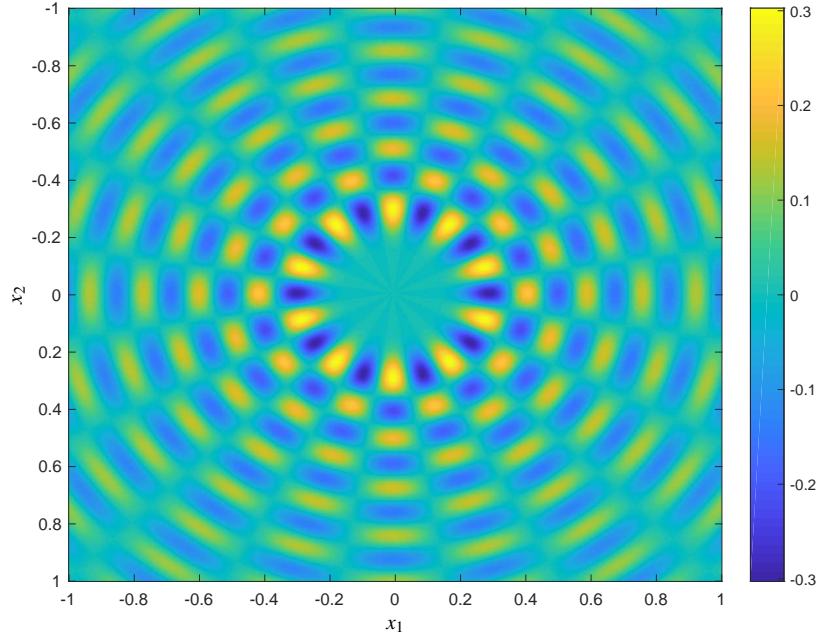} 
\caption{The real component of the regular wavefunction $\psi_{10}^i$ for $k=40$.}
\label{fig:inc_waves}
\end{figure}

We now compute numerical results for the particular example of a \emph{regular wavefunction} (synonymously referred to as \emph{Fourier-Bessel} functions), defined as
\begin{equation}\label{RWF}
\uinc (\bfx)=\psi_\ell^i(\bfx):=J_{|\ell|}(k|\bfx|)\e^{\imag\ell\theta_{\bfx}},\quad \ell \in \Z,
\end{equation}
where $J_\ell$ corresponds to the $\ell$th Bessel function of the first kind and $\theta_{\bfx}$ is the angle that $\bfx$ makes with the $x_1$-axis. The function $\psi_{10}^i$ for $k=40$ is depicted in Figure~\ref{fig:inc_waves}.

\begin{figure}[hbpt]
\centering
\begin{tikzpicture}\draw (-0.500000,-0.866025) -- (0.500000,-0.866025) -- (1.000000,-0.000000) -- (0.500000,0.866025) -- (-0.500000,0.866025) -- (-1.000000,0.000000) -- (-0.500000,-0.866025);\draw[red,thick,dashed,->] (-2.300000,0.000000) -- (-1.700000,0.000000);\draw[red,thick,dashed,->] (-1.150000,-1.991858) -- (-0.850000,-1.472243);\draw[red,thick,dashed,->] (1.150000,-1.991858) -- (0.850000,-1.472243);\draw[red,thick,dashed,->] (2.300000,0.000000) -- (1.700000,0.000000);\draw[red,thick,dashed,->] (1.150000,1.991858) -- (0.850000,1.472243);\draw[red,thick,dashed,->] (-1.150000,1.991858) -- (-0.850000,1.472243);\draw[thick,->] (-1.943937,1.229271) -- (-1.436823,0.908592);\draw[thick,->] (-2.292542,0.185073) -- (-1.694487,0.136793);\draw[thick,->] (-2.115953,-0.901523) -- (-1.563965,-0.666343);\draw[thick,->] (-1.454624,-1.781591) -- (-1.075157,-1.316828);\draw[thick,->] (-0.460059,-2.253518) -- (-0.340044,-1.665644);\draw[thick,->] (0.639900,-2.209192) -- (0.472970,-1.632881);\draw[thick,->] (1.593266,-1.658766) -- (1.177631,-1.226044);\draw[thick,->] (2.181634,-0.728336) -- (1.612512,-0.538336);\draw[thick,->] (2.270216,0.368946) -- (1.677985,0.272699);\draw[thick,->] (1.838718,1.381707) -- (1.359053,1.021262);\draw[thick,->] (0.985993,2.077936) -- (0.728777,1.535866);\draw[thick,->] (-0.092612,2.298135) -- (-0.068452,1.698621);\draw (0.000000,-1.241025) node{$\partial\Omega$};\end{tikzpicture}
\caption{Schematic diagram of Hexagonal scatterer (side length$=1$) used for Section {\ref{sec:radwaveEMB}}. 
For the Hexagon, it follows that $p=3$ and $q=4$, hence $\McansChSix=\NsidesChSix(q-1)=18$. 
The incident angles $A_{18}$ of canonical plane waves are shown; those chosen to coincide with $\Theta_*$ are {\color{red}dashed}.}
\label{HexDiag}
\end{figure}
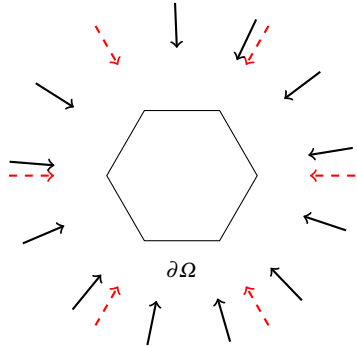

The Herglotz kernel (of Definition \ref{hergs}) of the $\ell$th regular wavefunction $\psi^i_\ell$, which follows by the Jacobi-Anger expansion (e.g. \cite[(3.89)]{CoKr:13}) and the definition \eqref{RWF}, is given by
\begin{equation}\label{FBherg1}
g_\ell(\alpha):=\left\{\begin{array}{cc}
(-1)^{\ell}\e^{-\imag \ell \alpha}/(2\pi\imag^\ell),&\ell\geq0,\\
\e^{-\imag \ell \alpha}/(2\pi\imag^\ell),&\ell<0,
\end{array}\right.
\end{equation}
hence we have $\psi_\ell^i=\uiHerg(\cdot;g_\ell)$. 
We now consider the case in which $\Omega$ is a regular hexagon {with unit side length (as depicted in Figure \ref{HexDiag}),} with wavenumber $k=1$. 
In our approximation we choose $\Nquad=20\times\max(k,\ell)$ equally spaced quadrature points {over the interval $[0,2\pi)$}, in order to resolve the oscillations of the integrand of \eqref{HhergDef}. 
These {quadrature points} are positioned such that
\begin{equation}\label{emb_quad_condish}
\min_{\underset{i\in\{1,\ldots,{\Nquad}\}}{\alpha_*\in\Theta_*}}|\tilde\alpha_i-\alpha_*|
\end{equation}
is maximised,
to avoid points where the error in the far-field approximation is significantly amplified (as observed in the peaks of Figure \ref{fixed_fig} for $\alpha=0$, and quantified by Lemma \ref{HatBound}).

Whilst the {solver's} DOFs per side of the scatterer are the same as for the smallest errors observed in Figure~\ref{fixed_fig}, for $\alpha=0$, in Figure~\ref{fig:HergGood} 
we observe global errors lower than were seen in the peaks of Figure~\ref{fixed_fig}. This is most likely because of the careful choice \eqref{emb_quad_condish}. Figure~\ref{RadErrsNaive} demonstrates how a naive approach is less stable, with much larger errors; choosing the quadrature points to be far from points of $\Theta_*$, in accordance with~\eqref{emb_quad_condish} does not lead to any significant improvement. 
This suggests that the method we present here is not only a means of removing numerical instability at certain points, but is essential for this generalised implementation of embedding formulae, which can compute the far-field coefficient for Herglotz wave function incident fields. 

\newcommand{\FBerrYaxis}{$\left|\mathcal F_\infty\psi^i_\ell(\theta) - \big[\mathcal{E}_{\mathcal{P}}^\circledast 
	\mathcal F_\infty \uiHerg(\cdot;\gHerg)\big](\theta)\right|/\|\mathcal F_\infty\psi^i_\ell\|_{L^2(0,2\pi)}$}

\newcommand{\FBerrYaxisBad}{$\left|\mathcal F_\infty\psi^i_\ell(\theta) - \big[\mathcal{E}_{\mathcal{P}} 
	\mathcal F_\infty \uiHerg(\cdot;\gHerg)\big](\theta)\right|/\|\mathcal F_\infty\psi^i_\ell\|_{L^2(0,2\pi)}$}
\begin{figure}[hbpt]
	\centering
	{\psfrag{[x]}{$\theta$}\psfrag{[y]}[Bc][Bc]{\FBerrYaxis}
		\psfrag{[x1]}[Bc][Bc]{$0$}\psfrag{[x2]}[Bc][Bc]{$\pi/2$}\psfrag{[x3]}[Bc][Bc]{$\pi$}\psfrag{[x4]}[Bc][Bc]{$3\pi/2$}\psfrag{[x5]}[Bc][Bc]{$2\pi$}\psfrag{[t]}[Bc][Bc]{$\ell=0,1,2$}
		\includegraphics[width=.95\textwidth]{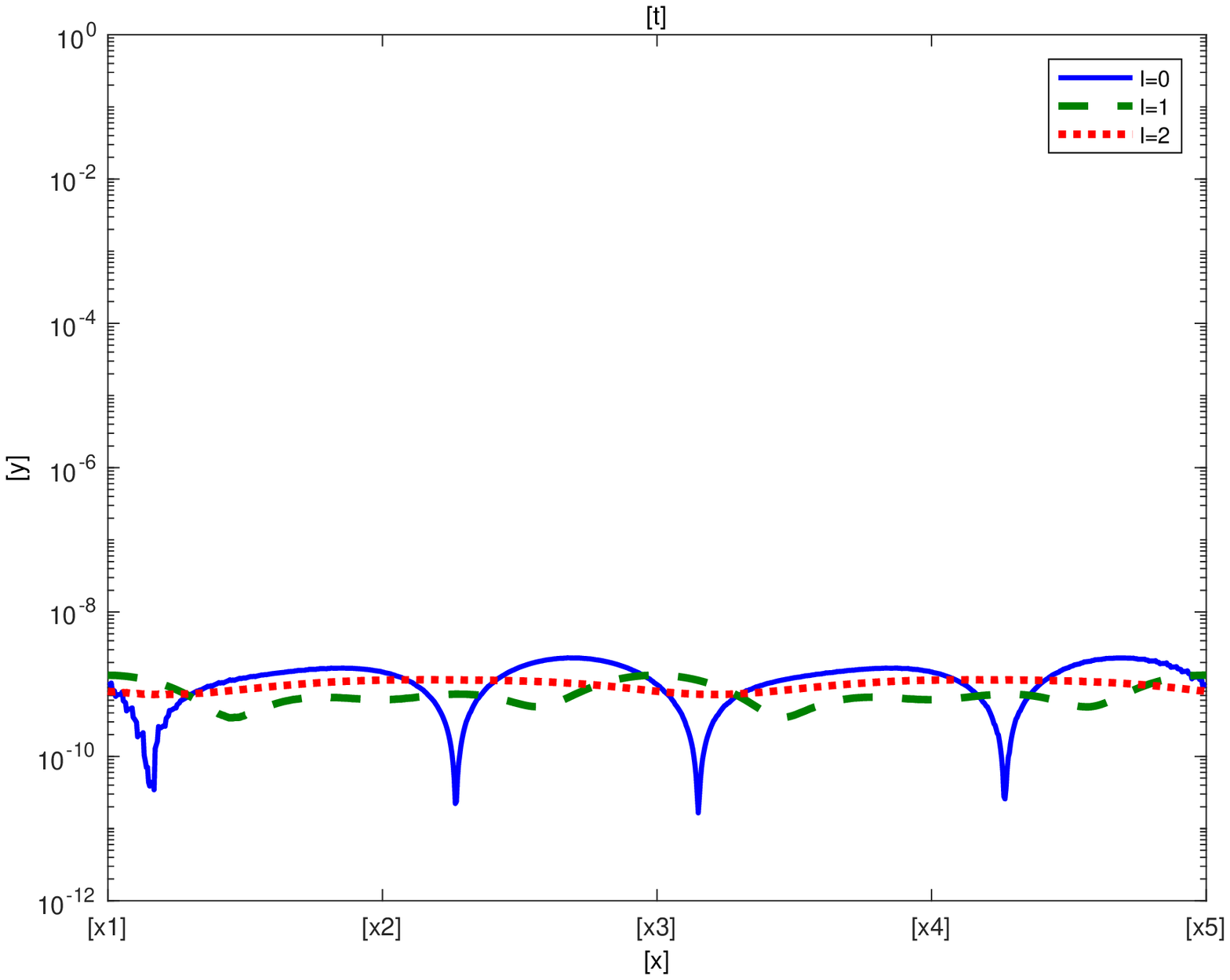}}
	{\psfrag{[x]}{$\theta$}\psfrag{[y]}[Bc][Bc]{\FBerrYaxis}
		\psfrag{[x1]}[Bc][Bc]{$0$}\psfrag{[x2]}[Bc][Bc]{$\pi/2$}\psfrag{[x3]}[Bc][Bc]{$\pi$}\psfrag{[x4]}[Bc][Bc]{$3\pi/2$}\psfrag{[x5]}[Bc][Bc]{$2\pi$}\psfrag{[t]}[Bc][Bc]{$\ell=3,4,5$}
		\includegraphics[width=.95\textwidth]{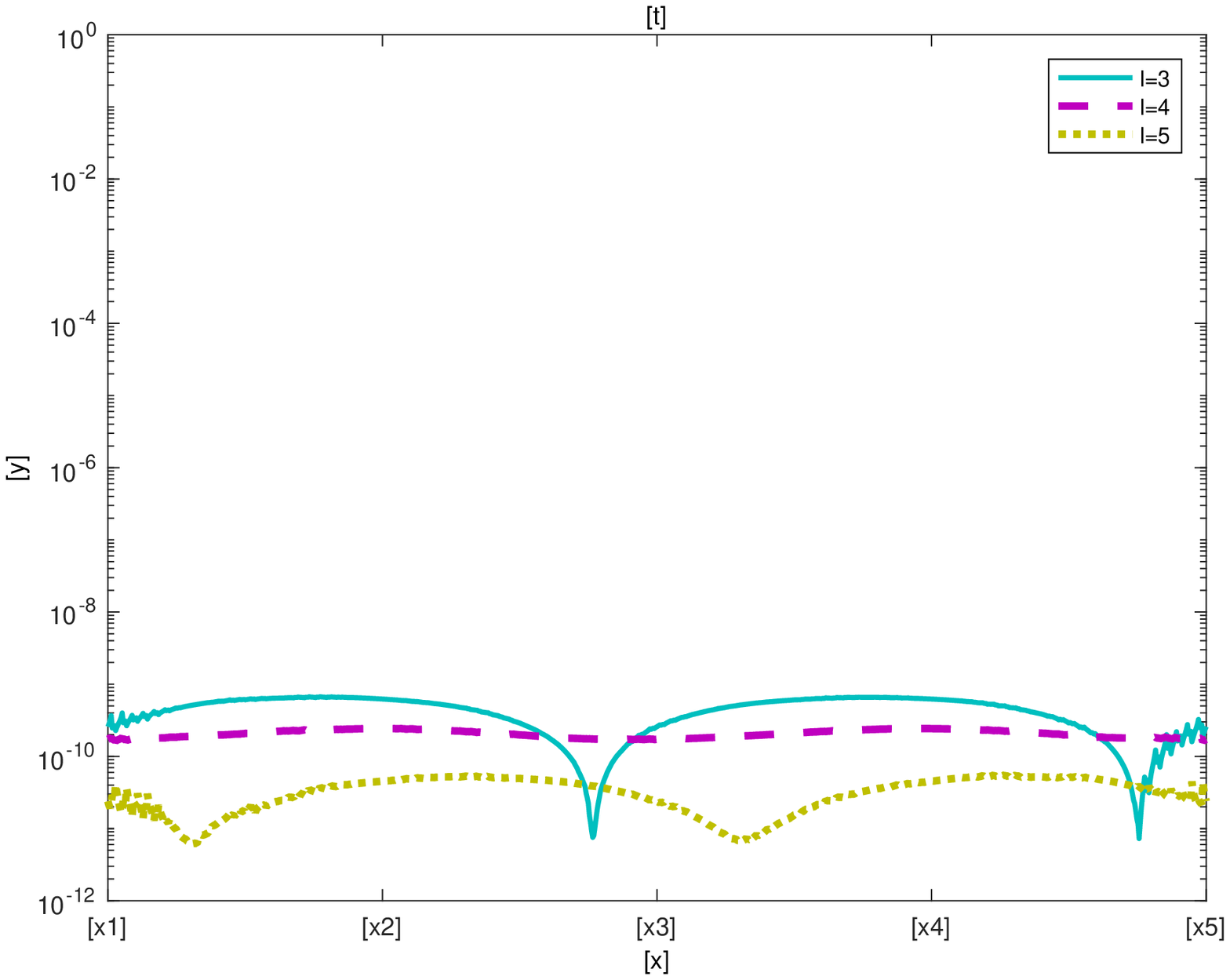}}
	\caption{ Error estimates for radiating wavefunctions \eqref{RWF} on the hexagon of Figure \ref{HexDiag}, using MPSpack with 90 DOFs as a reference solution, taking same DOFs per wavelength as for the square examples of Figure \ref{fig:heat}.}
	\label{fig:HergGood}
\end{figure}

\begin{figure}[hbpt]
	\centering
	{\psfrag{[x]}{$\theta$}\psfrag{[y]}[Bc][Bc]{\FBerrYaxisBad}
		\psfrag{[x1]}[Bc][Bc]{$0$}\psfrag{[x2]}[Bc][Bc]{$\pi/2$}\psfrag{[x3]}[Bc][Bc]{$\pi$}\psfrag{[x4]}[Bc][Bc]{$3\pi/2$}\psfrag{[x5]}[Bc][Bc]{$2\pi$}\psfrag{[t]}[Bc][Bc]{$\ell=0,1,2$}
		\includegraphics[width=\textwidth]{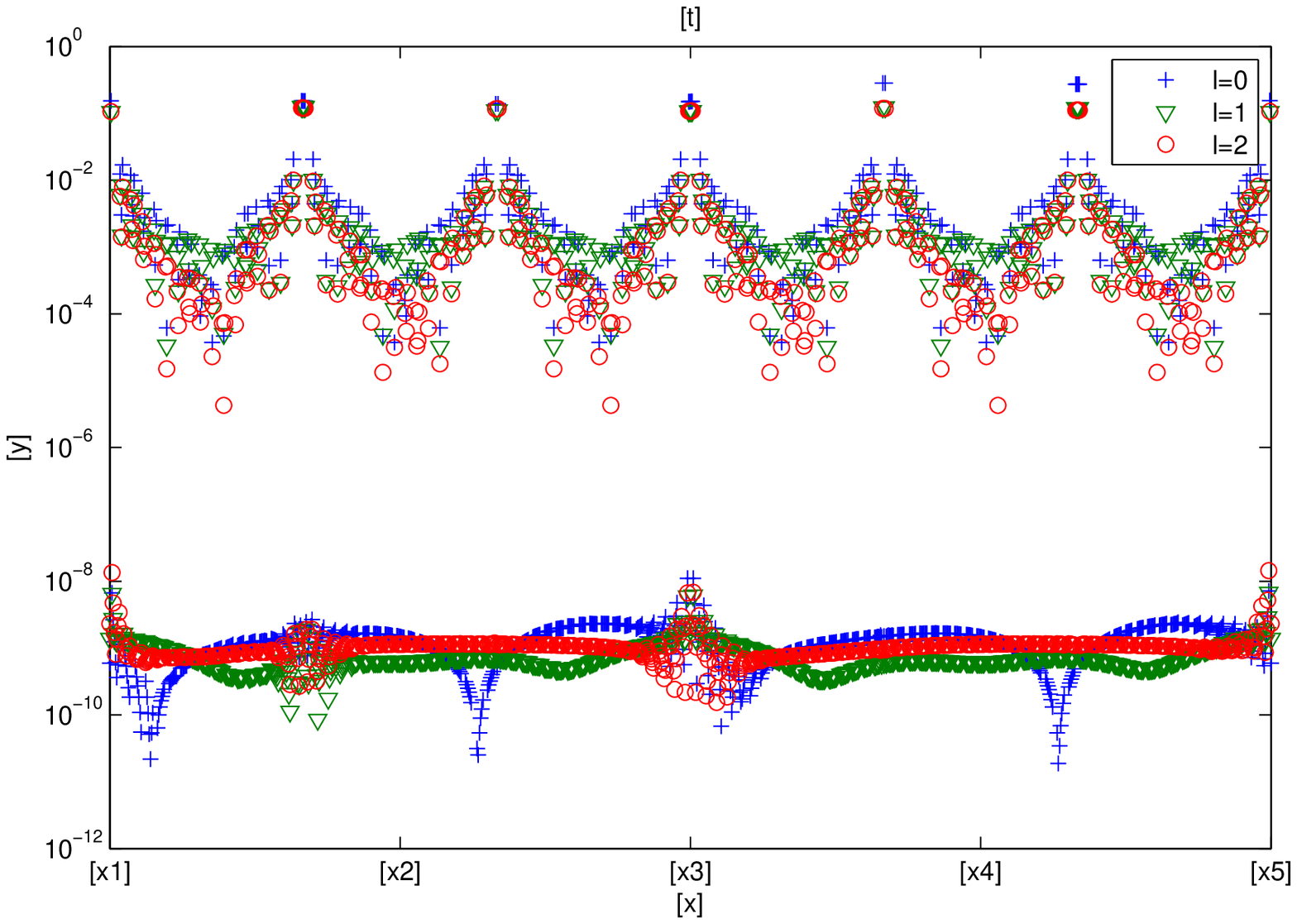}}
	{\psfrag{[x]}{$\theta$}\psfrag{[y]}[Bc][Bc]{\FBerrYaxisBad}
		\psfrag{[x1]}[Bc][Bc]{$0$}\psfrag{[x2]}[Bc][Bc]{$\pi/2$}\psfrag{[x3]}[Bc][Bc]{$\pi$}\psfrag{[x4]}[Bc][Bc]{$3\pi/2$}\psfrag{[x5]}[Bc][Bc]{$2\pi$}\psfrag{[t]}[Bc][Bc]{$\ell=3,4,5$}
		\includegraphics[width=\textwidth]{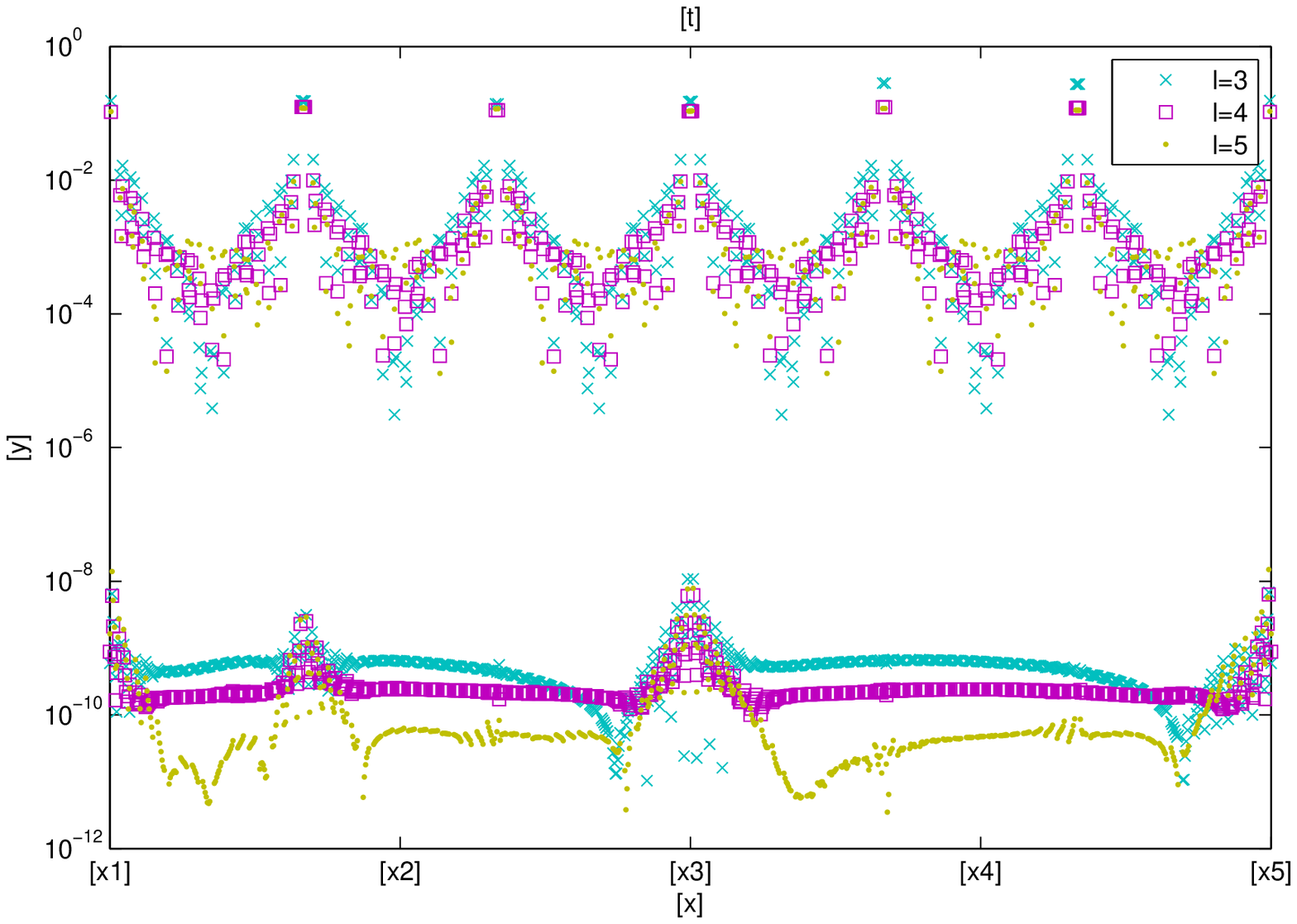}}
	\caption{  Error estimates for radiating wavefunctions \eqref{RWF} on the hexagon of Figure \ref{HexDiag}, using a completely naive embedding approach \eqref{emb0}, with a random selection of $A_\McansChSix$, standard midpoint rule used for weights and nodes for the integral, and no Taylor expansion.}\label{RadErrsNaive}
\end{figure}

This example demonstrates how, once the canonical problems have been solved, and certain derivatives of the far-field solutions are obtained, the far-field coefficient of any {Herglotz wave function incidence} can be computed relatively easily. Moreover, in certain cases, the global error can be reduced as part of this process.

\subsection{Application to T-matrix methods}

{T-matrix methods, first introduced in \cite{Wa:69}, are commonly used to solve scattering problems. They are particularly useful for moving body problems, and generalise easily from single-body to multiple-body problems. The method involves a decomposition of the incident field into a truncated basis of regular wavefunctions, and computation of a matrix $T$. This matrix maps a vector of regular wavefunction coefficients to a vector of coefficients corresponding to a decomposition of the scattered field in terms of a truncated basis of radiating wavefunctions (these are defined similarly to \eqref{RWF}, with the Bessel function $J_{|\ell|}$ replaced by a Hankel function of the same order).}

\begin{figure}[hbpt]
	\includegraphics[width=\textwidth]{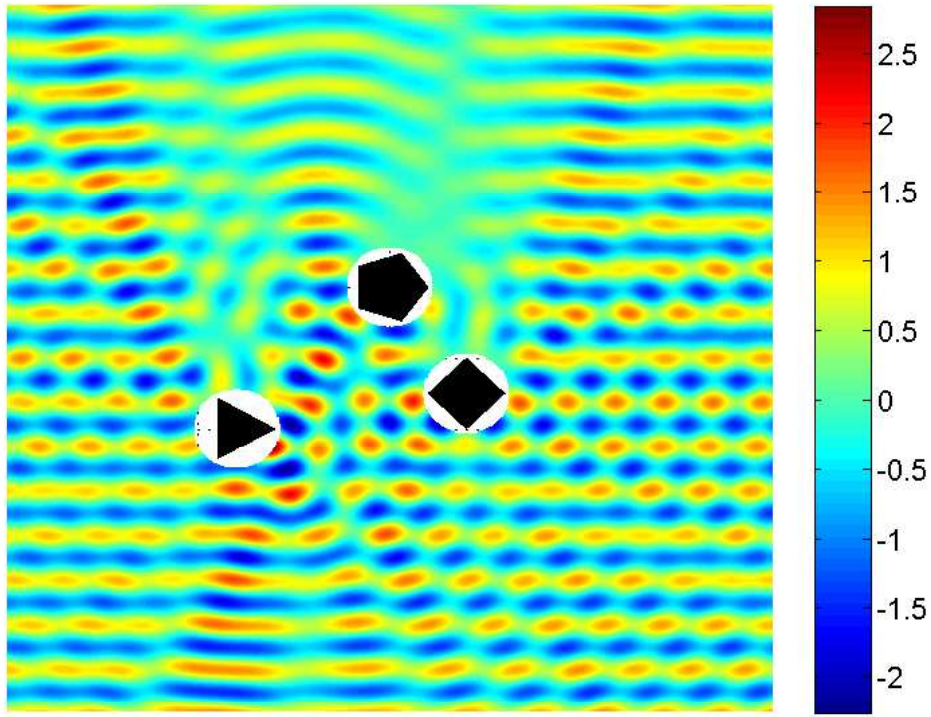}
	\caption{Real part of total field for a configuration of multiple polygons, with incident field ${u^i_{3\pi/2}}$ solved using Tmatrom coupled with MPSpack and the combined embedding approximation (of Definition \ref{s:CEA}) as the solver used for the embedding implementation, which in turn is used as the solver for Tmatrom. The representation is only valid outside of the union of balls containing each obstacle, a fundamental property of T-matrix methods.}\label{fig:TmatEG}
\end{figure}

{In 2017 Ganesh and Hawkins released} \emph{Tmatrom}, a numerically stable T-matrix software package \cite{GaHa:16}{ to compute $T$}. This method deviates from standard (but often unstable) approaches such as \cite{Wa:69}, in that it requires the user to input the far-field coefficient for regular wavefunctions $\{\psi_\ell\}_{\ell=-\MtruncChSeven}^{\MtruncChSeven}$, where
\begin{equation}\label{TmatromCondish}
\MtruncChSeven =\left\lceil kR_-+4(kR_-)^{1/3}+5\right\rceil
\end{equation}
{and $R_-$ is the radius of a ball containing $\Omega$ (three such balls can be found in Figure \ref{fig:TmatEG}).}
The condition \eqref{TmatromCondish} guarantees stability of the Tmatrom method \cite{GaHaHi:12}. Moreover, accuracy of the Tmatrom method is maintained as the wavenumber $k$ increases, provided that \eqref{TmatromCondish} holds. Hence the software requires (an approximation of) the far-field coefficient induced by $O(k)$ distinct regular wavefunctions, namely
\begin{equation}\label{Treq}
\{\mathcal{F}_\infty\psi_\ell^i\}_{\ell=-\MtruncChSeven}^{\MtruncChSeven}.
\end{equation}
Rather than solve $2\MtruncChSeven+1=O(k)$ problems of regular wavefunction incidence, we will now show how embedding formulae can be used to reduce this to a number (of plane wave problems) which does not grow with wavenumber $k$. By combining \eqref{emb0}, \eqref{HhergDef} and \eqref{FBherg1} we can write for $\ell=-\MtruncChSeven,\ldots,\MtruncChSeven$,
\begin{align}\label{OktoO1}
\mathcal{F}_\infty\psi_\ell^i(\theta)&=\int_0^{2\pi}g_\ell(\alpha)\frac{\sum_{m=1}^\McansChSix B_m(\alpha)\hat{D}(\theta,\alpha_m)}{\Lambda(\theta,\alpha)}\dd{\alpha}\nonumber\\
&=\int_0^{2\pi}g_\ell(\alpha)\frac{\sum_{m=1}^\McansChSix B_m(\alpha)\Lambda(\theta,\alpha_m){\mathcal{F}_\infty u^i_{\alpha_m}(\theta)}}{\Lambda(\theta,\alpha)}\dd{\alpha}.
\end{align}
Hence, we may solve for $\McansChSix=O(1)$ incident plane waves to obtain (by inserting into \eqref{OktoO1}) an approximation to the far-field coefficient for any $\ell$. Then any far-field pattern of the set \eqref{Treq} follows immediately from the representation \eqref{OktoO1}, no further solves are required, one can simply re-integrate \eqref{OktoO1} against a different Herglotz kernel $g_\ell$. Due to the numerical instabilities at certain observation angles, in practice one must use the approximation \eqref{HhergDef}-\eqref{HapproxDef} with $\gHerg=g_\ell$, as discussed in \S\ref{sec:radwaveEMB}, hence the required inputs to Tmatrom, which approximate the elements of \eqref{Treq}, are
\[
\left\{[\mathcal{E}^\circledast_{\mathcal{P}}\mathcal{F}_\infty \uiHerg (\cdot;g_\ell)]\left(\theta,\gHerg; \MTayChSix, ({w}_i)_{i=1}^{\Nquad},(\tilde\alpha_i)_{i=1}^{\Nquad}\right)\right\}_{\MtruncChSeven=-\ell}^\MtruncChSeven,
\]
which can be obtained by solving only $\McansChSix$ problems of plane wave incidence.
Figure \ref{fig:TmatEG} shows the output of the combination of Tmatrom with our embedding solver and the MPSpack solver, for the problem of plane wave scattering with three regular polygons. This required $8$, $12$ and $30$ solves on the triangle, square and pentagon respectively, a total of $50$ solves {for a single scatterer each}, a number independent of $k$. For wavenumber $k=5$, using MPSpack without solving via the embedding formulae results in a total of $27$ solves on each scatterer, hence a total of $81$ solves. So even at a relatively low wavenumber, the embedding formulae can reduce the number of solves required by Tmatrom; 
as $\MtruncChSeven$ grows with $O(k)$, the number of solves required by embedding {formulae} is independent of $k$.  There is one hidden cost with embedding {formulae}, namely that the number of terms required in the Taylor series (\ref{embPT}) and (\ref{shiTaylor}) in order to maintain accuracy will need to grow with $k$; we leave detailed consideration of this to future work.

\begin{acknowledgements}
The authors would like to thank Nick Biggs and Stuart Hawkins for many helpful discussions.
\end{acknowledgements}

\bibliographystyle{apa}
\bibliography{AllRefsV3}   

\end{document}